\documentclass[11pt,oneside]{amsbook}
\usepackage{amsmath}
\usepackage{breqn}
\usepackage{pdfpages}
\usepackage{multicol}
\usepackage{nicematrix,tikz}
\usepackage{pgfplots}
\pgfplotsset{compat=1.8}
\usepgfplotslibrary{fillbetween}
\usepackage{collectbox}
\usepackage{cancel}
\usepackage{float}
\usepackage[none]{hyphenat}

\usepackage{array}
\usepackage{arydshln}
\usepackage[hidelinks]{hyperref}

\usepackage{amssymb}

\usepackage{skak}
\usepackage{graphicx}
\usepackage{tikz}

\usepackage{mathtools}

\DeclarePairedDelimiter\floor{\lfloor}{\rfloor}

\usepackage{pdflscape}
\usepackage{setspace}

\usepackage{stackengine}
\newlength\Lobj
\newsavebox{\newobj}

\newcommand\LoF[2]{%
  \def\obj{$#1$~$#2$\,}%
  \setlength\Lobj{\widthof{\obj}}%
  \sbox\newobj{\stackon[3pt]{\obj}{\rule{\Lobj}{.3pt}}}%
  \usebox{\newobj}%
  \rule{.3pt}{\ht\newobj}%
  }

\numberwithin{equation}{chapter}
\numberwithin{figure}{chapter}

\theoremstyle{plain} 

\begin{document}
\frontmatter 


\centerline{\huge \textbf{Primes Between Squares}:}
\medskip
\centerline{\LARGE \textbf{Commentary on Appendix 8 of \textit{Laws Of Form}}}

\bigskip

\begin{center}
 
 \LoF{\vphantom{@}}{\hphantom{.}} \LoF{\vphantom{@}}{\hphantom{.}} = \LoF{\vphantom{@}}{\hphantom{.}}\\
 
 \medskip
 
 \LoF{\LoF{\vphantom{@}}{\hphantom{.}}}{\!}  =  \hphantom{\LoF{\vphantom{@}}{\hphantom{.}} \LoF{\vphantom{@}}{\hphantom{.}}} \\

\end{center}

\centerline{\Huge }
\medskip
\bigskip

\centerline{\LARGE J. M. Flagg}
\medskip
\medskip
\medskip
\centerline{\LARGE Louis H. Kauffman}
\begin{center}
Department of Mathematics, Statistics, and Computer Science,\\
851 South Morgan Street,\\
University of Illinois at Chicago,\\
Chicago, IL 60607-7045\\
and\\
International Institute for Sustainability wih Knotted Chiral Meta Matter,\\
WPI-SKCM2, Hiroshima University,\\ 
1-3-1 Kagamiyama, Higashi-Hiroshima, Hiroshima 739-8526, Japan\\
\textbf{loukau@gmail.com}
\end{center}
\medskip
\medskip
\centerline{\LARGE Divyamaan Sahoo}
\begin{center}
Department of Electrical \& Computer Engineering,\\
University of Massachusetts Dartmouth,\\
North Dartmouth, MA 02747\\
and\\
Graduate Program in Acoustics,\\
The Pennsylvania State University,\\
University Park, PA 16802\\
\textbf{divyamaansahoo@gmail.com}

\end{center}



\newpage
\thispagestyle{empty}

\vspace*{0.5 cm}

\begin{center}
    \begin{quote}
    \noindent \textit{An argument follows.}\\
    \medskip
    \noindent \textbf{James Joyce, \textit{Finnegans Wake}}
\end{quote}

\bigskip
\bigskip

\begin{quote}
    \noindent \textit{Hail Thoth, architect of truth, give me words of power that I may form the characters of my own evolution. I stand before the masters who witnessed the genesis, who were the authors of their own forms, who rolled into being, who walked the dark circuitous passages of their own becoming, ...}
\end{quote}

\begin{quote}
    \noindent \textit{Hail Thoth, architect of truth, give me words of power that I may intuit the symbols of dream and command my own becoming. I stand before the masters who witnessed the working of magic, ...}\\
\end{quote}

\begin{quote}
    \noindent \textit{Hail Thoth, architect of truth, give me words of power that I may complete my story and begin life anew. I stand before the masters who witnessed the plowing of earth, who saw the seed that entered the fields spring into corn and barley, who sent the flood and sun, who saw men among the wheat swinging scythes, who saw women baking bread. ...}\\
\end{quote}

\begin{quote}
    \noindent \textit{Hail Thoth, architect of truth, give me words of power that I may tell the truth of my own becoming. I stand before the masters who witness the judgment of souls, who sniff out the misdeeds, the imperfections, the lies and half-truths we tell ourselves in the dark.} ...

\medskip
\medskip

    \noindent \textbf{Normandi Ellis, \textit{Awakening Osiris}}
\end{quote}
\end{center}

\newpage

\mainmatter 

\chapter*{Primes Between Squares:\\ Commentary on Appendix 8 of \textit{Laws Of Form}}



\noindent \textbf{Abstract:}\\
There are two conjectures regarding the distribution of primes between squares attributed to Legendre, a weak and a strong conjecture. The weak conjecture\footnote{When reformulated, this conjecture is the fourth problem of Edmund Landau, who in his lecture ``Gel\"{o}ste und ungel\"{o}ste Probleme aus der Theorie der Primzahlverteilung und der Riemannschen Zetafunktion" at the International Congress of Mathematicians in Cambridge, England, on August 23\textsuperscript{rd}, 1912, characterised his four problems as ``unattackable at the present state of science". In Landau's words, ``Liegt zwischen $n^2$ und $(n+1)^2$ f\"{u}r alle positiven ganzen $n$ mindestens eine Primzahl?" (Is there at least one prime number in between $n^2$ and $(n+1)^2$ for all positive whole $n$?)} is that there is \textit{at least one} prime between consecutive squares. The strong conjecture is that there exist \textit{at least two} primes between consecutive squares. In Appendix 8 of \textit{Laws Of Form} (Bohmeier Verlag, revised sixth English edition, 2015), Spencer-Brown sets out to prove the stronger of the two. In unpublished material\footnote{The authors are indebted to Andrew Crompton, Graham Ellsbury, Thomas Wolf, and Moshe Klein for insights and conversations.}, he provides stronger and sharper theorems and conjectures on the regular and predictable increase of the number of primes between consecutive squares. In Appendix 7 of \textit{Laws Of Form} (Bohmeier Verlag, 2015), and in \textit{Gesetze der Form} (Bohmeier Verlag, 1997), he claims:

\fbox{\begin{minipage}{32em}

\noindent The total number $t(n)$ of primes between $n^2$ and $(n+1)^2$, where $n$ is a natural number greater than $1$, is within the limits $$A-(B-1) < t(n) < A+(B-1),$$ where $A$ is $\frac{n}{\log n}$ and $B$ is $\frac{A}{\log A}$.

\end{minipage}}



\noindent In other words, the number of primes between squares is asymptotic to $\frac{n}{\log n}$. This article provides commentary and context to the number theory of Spencer-Brown.


\setcounter{section}{-1}

\newpage
\thispagestyle{empty}
\section{Introduction}

\noindent This paper provides a commentary and guide to Appendix 8 of \textit{Laws Of Form}, which is an article by Spencer-Brown on number theory and his proofs of the conjecture that there are at least two prime numbers between any consecutive squared numbers $n^2$ and $(n+1)^2$. For example, $1^2=1$ and $2^2 = 4$ and the numbers between them are $2$ and $3$ thus we have two prime numbers between $1^2$ and $2^2$. This phenomenon continues for all pairs $n^2$ and $(n+1)^2$ and was a conjecture in number theory since Legendre. In Spencer-Brown’s appendix he gives his proofs of the conjecture. We shall see that those proofs are a highly original mixture of standard rigorous arguments and also some stated facts about the way numbers behave that would be considered conjectures by most number theorists. These phenomena are very interesting and constitute a deep observation about the nature of number itself. We hope that our guide will enable the reader to gain insight into Spencer-Brown’s point of view.

\noindent \textbf{Section 1} defines the square segment and introduces the stroke notation of counting primes in segments. Here, we introduce the sharp conjectures made by Spencer-Brown on the distribution of primes between squares, beginning with Theorem 1, which is a \linebreak reformulation of Legendre's weak conjecture that there exists one prime between \linebreak squares. The fundamental observation is that the primes $< n$ determine the primes $< n^2$, and this guides Lemmas 1 and 2, which are both proved rigorously.

\noindent \textbf{Section 2} lays the foundation of Spencer-Brown's approach to the Sieve of Eratosthenes by clarifying his use of the word ``strike" as an operator and operation, culminating in the ``prime paradigm", which is a repeating pattern of composites generated by a set of primes. This pattern has a four-fold symmetry about multiples of the primorial (product of primes) and multiples of half the primorial.

\noindent \textbf{Section 3} defines a sect as a sequence of consecutive composites and reformulates Legendre's conjecture (Theorem 1) by asking the question: what is the maximum sect that can appear in a segment? Two conjectural lemmas that are closely related to Theorem 1 are introduced here, namely Lemma 3, which asserts that maximal odd sects assume one of two forms, quadrantic or reflexive, and Lemma 4, which asserts that the quadrantic sect length is always lesser than or equal to the reflexive sect length. This section ends with a rigorous proof of Theorem 1, however, this proof rests on the validity of Lemmas 3 and 4.

\noindent \textbf{Section 4} generalizes Theorem 1 to a powerful family of conjectures of the following form: the next prime after $p$ is at most $\frac{1}{k}\sqrt{p}$ away, where $k=\frac{1}{2}$ and for primes $p=7, 113, 1327, ...$, the value of $k$ jumps to a larger value from which it cannot slip back. This is called the ``ratchet" phenomenon. The section ends with two of Spencer-Brown's stronger conjectures: a Corollary, which states that for any $n>131$ there is at least one prime in $(n-\sqrt{n},n)$ and $(n,n+\sqrt{n})$, and Theorem 2, which is a reformulation of Legendre's strong conjecture that there exist at least two primes between squares.

\noindent \textbf{Section 5} introduces Spencer-Brown's theory of elementals and connects it with his theory of modulators, rejoining his algorithm in Lemma 3 to the initials of \textit{Laws Of Form}.

\noindent Finally, the paper ends with three appendices. \textbf{Appendix 1} provides unpublished excerpts of Spencer-Brown on sharp bounds and estimates to $\pi(x)$. \textbf{Appendix 2} contextualizes Spencer-Brown's conjectures in a set of related conjectures and theorems on the distribution of prime numbers. \textbf{Appendix 3} illustrates the sole exception to Spencer-Brown's Lemma 3 at $23$. 

\noindent Following Petzold's convention in \textit{The Annotated Turing} (2008), this article uses a box when quoting excerpts from Spencer-Brown in Appendix 8 of \textit{Laws Of Form} (Bohmeier Verlag, 2015) accompanied by elucidatory annotations and corrections.

\newpage
\thispagestyle{empty}
\section{Preliminaries}


\noindent Appendix 8 ``Primes between squares", \textit{Laws Of Form}, begins with the definition of ``square segment" or ``segment" as an interval between squares whose roots differ by $1$. Hence, $\text{seg}(x)$ indicates the stretch of natural numbers between $x^2$ and $(x+1)^2$. Call the square segment $\text{seg}(n)$ when the roots are natural numbers. For example:

\vspace{-0.4 cm}

\begin{center}
$\text{seg}(2) = \{5,6,7,8\}$,\\
$\text{seg}(\sqrt{2}) = \{3,4,5\}$.
\end{center}

\noindent Note that $\text{seg}(2)$ is the set of natural numbers between $2^2$ and $3^2$, i.e., between $4$ and $9$, and $\text{seg}(\sqrt{2})$ is the set of natural numbers between $(\sqrt{2})^2$ and $(\sqrt{2}+1)^2$, i.e., between $2$ and $\approx 5.8284$. Note, $\text{seg}(0) = \{ \}$, the stretch of natural numbers between $0^2$ and $1^2$. 

\noindent Adapting Landau's notation $\pi(x)$ for the number of primes $\leq x$, Spencer-Brown defines $\hphantom{.}^\circ\pi\text{seg}(x)$ read ``stroke pi seg of $x$"\footnote{Note: the $\hphantom{.}^\circ$ or `stroke' notation is used before by Spencer-Brown in the context of modulators (reductors), Chapter 11, \textit{Laws Of Form}, page 66 - ``Let a marker be represented by a vertical stroke, thus $\vert$" - and in Spencer-Brown's unpublished manuscript, ``Introduction to Reductors", Section 1.6: Conventions, page 10 - ``The standard \textit{drive} D to a four-stroke reductor is taken to be 1010, and each unit is called a \textit{stroke}, written $^\circ 1, ^\circ 2, ^\circ 3$, or $^\circ 4$". So far, Spencer-Brown's earliest known use of the term ``stroke" is in ``The Universal Operators in Logic". In this unpublished manuscript, he examinines other primitive mathematical operators, of Sheffer, Stamm, Peirce, Nicod, et alia, all of which he deems inadequate due to restrictions and constraints delimiting scope of operations. He raises the question why the bar (a bar drawn over an algebraic equation) had not been considered as a ``sole" operator and writes that the theoretical significance of its practical form and applicability were missed. The bar just needed some terminal modification so that it could be read correctly. He says, ``To overcome the mechanical difficulty, it is suggested that the bar should be determined by a stroke, thus, \LoF{\vphantom{@}}{\hphantom{.}}". It is also interesting to note that the Theory of Types, whether simple or ramified, uses right upper superscripts to denote a proposition’s type and place in an elevated hierarchy to elide and avoid paradox, thus eliminating  any formal system incapable of referring to itself. Spencer-Brown possibly adapted this notation by moving the superscript to the left of the operator/operand and designating it as a circle (or ouroboros) of modular degree with valid closure, referring to and calculating its contents recursively, and calling it a ``stroke."} or, equally often, ``pi seg of $x$" to denote the number of primes in $\text{seg}(x)$. So, $\pi(x)$ goes to $^\circ \pi(x)$. At the outset, Spencer-Brown claims:\\
\noindent \textbf{i.} When $x$ is a real number $>$ the least absolute value of $(\sqrt{2}-1), \hphantom{.}^{\circ}\pi\text{seg}(x) \geq 1$. \\
\textbf{ii.} When $x$ is a natural number $n$, $^{\circ}\pi\text{seg}(n) \geq 2$. \\ \textbf{iii.} When $n \geq 6, \hphantom{.} ^{\circ}\pi\text{seg}(n) \geq 3$, and the minimum number of such primes, or more exactly its lesser bound, will increase regularly and predictably as $n$ increases. 

\noindent In this context, Spencer-Brown is referring to his Prime Limit Theorem, Equation (5) on pg. 182, Appendix 7, \textit{Laws Of Form} (Bohmeier Verlag, revised sixth English edition, 2015), that $\hphantom{.}^\circ\pi\text{seg}(n)$ is asymptotic to $\frac{n}{\log n}$.

\noindent After defining the (square) segment, Spencer-Brown introduces Theorem 1, which is a reformulation of Legendre's weak conjecture.

\fbox{\begin{minipage}{32em}
\noindent \textbf{Theorem 1.}\\
There is at least one odd prime in the square segment of every natural number.
\end{minipage}}

\noindent Spencer-Brown proves Theorem 1 by introducing the following Lemmas 1 through 4. \textit{Note: Lemma 1A is not provided by Spencer-Brown but is included here for clarity.}

\fbox{\begin{minipage}{32em}
\noindent \textbf{Lemma 1}:\\ 
A number $h$ is prime iff $h$ is an integer and $h$ does not make an integer quotient with any integer $i$ in case $2 \leq i \leq \sqrt{h}.\hphantom{.}$\footnote{The constraint on $i$ ensures that $i$, $\sqrt{h}$, and $h$ are all positive, but in fact primes, like moduli, are signless.}\\
\textbf{Proof:} There is no such $i$ unless $\sqrt{h} \geq 2$. If $i$ divides $h$ then $h/i = j$, say. Now either $i = j$, $i<j$, or $i>j$. In the first case $i = \sqrt{h}$, in the second $i < \sqrt{h}$, and in the third $j<\sqrt{h}$ and can replace $i$ as the divisor. So effectively all proper divisors of $h$ will have been tried by confining $i$ in case $2 \leq i \leq \sqrt{h}$.\qed
\end{minipage}}

\noindent In the footnote to Lemma 1, Spencer-Brown is referring to the modular number system. This powerful Lemma states that in order to test whether any number is a prime, it is sufficient to look for factors only up to the square-root of that number. This fact is well known.

\noindent \textit{NOTE: The following lemma is not provided by Spencer-Brown.}\\
\textbf{Lemma 1A}:\\
The least prime factor of any odd composite between $n^2$ and $(n+1)^2$ is $\leq n$.\\
\textbf{Proof:} The primes $p_1, p_2, p_3, \ldots, p_j$ solely determine all the primes up to $p_{j+1}^2$ (see Volume I, Chapter XIII, \textit{History of the Theory of Numbers}, L. E. Dickson). If $n = p_j$, then $(n+1)^2 < p_{j+1}^2$. If $(n+1) = p_{j+1}$, then $(n+1)^2 \leq p_{j+1}^2$. In both cases, the primes $p_1, p_2, p_3, \ldots, p_j \leq n$ determine all the primes up to $(n+1)^2$, so the least prime factor of any odd composite between $n^2$ and $(n+1)^2$ is $\leq n$. \qed


\fbox{\begin{minipage}{32em}
\noindent \textbf{Lemma 2}:\\
There are $n$ odd numbers between $n^2$ and $(n+1)^2$.
\end{minipage}}

\noindent \textbf{Proof:} $(n+1)^2 - n^2 = 2n+1$. In a stretch of $2n+1$ integers, $n$ must be odd.\qed

\noindent First consider, seg$(1) = \{2,3\}$, seg$(2) = \{5,6,7,8\}$, seg$(3) = \{10,11,12,13,14,15\}$, and seg$(4)=\{17,18,19,20,21,22,23,24\}$. Carefully observe the odd numbers in seg$(5)$:

\begin{center} \begin{tabular}{c c c : c c c c}
        \LARGE $5^2$ & $27$ & $29$ & $31$ & $33$ & $35$ & \LARGE $6^2$\\
        $\square$ & $\bullet$ & $\bullet$ & $\bullet$ & $\bullet$ & $\bullet$ & $\square$\\
        $5$ & $3$ & & & $3$ & $5$ &
    \end{tabular}
\end{center}



\noindent The dashed line marks the primorial $2\times 3\times 5$. The numbers $5$ and $3$ in the third line indicate the first prime that divides the odd number above it, so $29$ and $31$ are prime.\\
\noindent By Lemma 2, there are $5$ odd numbers between $5^2$ and $6^2$.\\
\noindent By Lemma 1A, the least prime factor of any odd composite in seg$(5)$ is $ \leq 5$.\\
\noindent By Lemma 1, $h=29$ is prime because it does not make an integer quotient $i$ with any integer $i$ in case $2 \leq i \leq \sqrt{29}\approx 5.385$, i.e., $i \in \{2, 3, 4, 5\}$. Similarly, $h=31$ is prime.

\noindent Call all primes $\leq n$ the \textit{prime generators of $n^2$} or \textit{pg of $n^2$} for short. Alternatively, the primes $\leq \sqrt{n}$ are the \textit{prime generators} of $n$. This idea anticipates the recursive relationship of the prime counting function\footnote{In a footnote on page 190, Spencer-Brown writes, ``The true pg of $n$ is \textit{all} the multiplicative integers $\leq \sqrt{n}$, because it demands no previous knowledge of primes". The Legendre/Spencer-Brown formula in Appendix 9 of \textit{Laws Of Form} shows how the prime counting function is self-referential. The formula is $\pi (n) = \pi ( \sqrt{n}) + \sum_{d = 1}^{\sqrt{n}} \mu (d) \floor{\frac{n}{d}}$, where $d \in\{$the set of all combinations of prime generators $\leq \sqrt{n}$, including $1\}$. For $d=1$, Spencer-Brown changes $\floor{\frac{n}{d}}$ to $\floor{\frac{n-1}{d}}$. See ``Laws Of Form and the Riemann Hypothesis", \textit{Laws Of Form: A Fiftieth Anniversary}, World Scientific, 2021.}, specifically how $\pi(n)$ is a function of $\pi(\sqrt{n})$. This relationship is illustrated by Spencer-Brown in Appendix 9 of \textit{Laws Of Form} (Bohmeier Verlag, 2015).

\newpage
\thispagestyle{empty}
\section{The Foundation}

\noindent On page 190 of Appendix 8, Spencer-Brown lays the foundation to his approach:

\fbox{\begin{minipage}{32em}
Imagine the system of integers, call it an arithmetic $A$, mapped onto a line of regularly-spaced indistinguishable marks or points stretching endlessly in either direction. Call these marks the elements of a system $S$.
\end{minipage}}

\begin{center}
    \begin{tabular}{c c c c c c c c c c c c}
        $\ldots$ & $\bullet$ & $\bullet$ & $\bullet$ & $\bullet$ & $ \bullet $ & $\bullet$ & $\bullet$ & $\bullet$ & $ \bullet $ & $\bullet$ & $\ldots$
    \end{tabular}
\end{center}


\fbox{\begin{minipage}{32em}
\noindent Consider prime multisectors $d_1, d_2, d_3, \dots$ that can be imbedded in every second, every third, every fifth, and so on of the elements in $S$.
\end{minipage}}

\noindent Multisector $d_i$ is equal to $p_i$ (the $i$\textsuperscript{th} prime), i.e. $d_1 = 2, d_2 = 3, d_3 = 5, ...$. Hence, fix an origin and mark the prime $p_i$ (and all integer multiples of $p_i$) with $d_i$ as indicated:

\begin{center}
    \begin{tabular}{c c c c c c c c c c c c c}
        $\ldots$ & $\bullet$ & $\bullet$ & $\bullet$ & $\bullet$ & $\bullet$ & $ \bullet $ & $\bullet$ & $\bullet$ & $\bullet$ & $ \bullet $ & $\bullet$ & $\ldots$\\
         & 0 & & $d_1$ & & $d_1$ & & $d_1$ & & $d_1$ & & $d_1$ & 
    \end{tabular}
\end{center}

\begin{center}
    \begin{tabular}{c c c c c c c c c c c c c}
        $\ldots$ & $\bullet$ & $\bullet$ & $\bullet$ & $\bullet$ & $\bullet$ & $ \bullet $ & $\bullet$ & $\bullet$ & $\bullet$ & $ \bullet $ & $\bullet$ & $\ldots$\\
         & 0 & & & $d_2$ & & & $d_2$ & & & $d_2$ & & 
    \end{tabular}
\end{center}

\begin{center}
    \begin{tabular}{c c c c c c c c c c c c c}
        $\ldots$ & $\bullet$ & $\bullet$ & $\bullet$ & $\bullet$ & $\bullet$ & $ \bullet $ & $\bullet$ & $\bullet$ & $\bullet$ & $ \bullet $ & $\bullet$ & $\ldots$\\
         & 0 & & & & & $d_3$ & & & & & $d_3$ & 
    \end{tabular}
\end{center}

\fbox{\begin{minipage}{32em}
... the multisectors are not primes or prime divisors, but merely operators on elements of $S$ that can be placed where we please to imbed with or strike certain of these elements.
\end{minipage}}

\noindent In the handwritten edition of \textit{Primes Between Squares}, Spencer-Brown uses the terms ``multisector" and ``decimator" interchangeably. He writes, ``Since decimation is a multisection by ten, we may conveniently extend the meaning of this more-familiar term and call each multisector a decimator, giving it the name of the prime divisor whose span it is associated with". The following two pages excerpt page 190 of Appendix 8.

\fbox{\begin{minipage}{32em}
\noindent Consider $S$ after $d_1$ has struck every second element, and eliminate the elements that have been struck, leaving the remainder to form a system $S'$, say. $S'$ is now identical to $S$, and $d_1$ is the only multisector that makes such an identity. So we can save labour by mapping the odd numbers in $A$, call them $A'$, onto the elements of $S'$.
\end{minipage}}

\begin{center}
    \begin{tabular}{c c c c c c c c c c c c c}
        $\ldots$ & $\bullet$ & $\bullet$ & $\bullet$ & $\bullet$ & $\bullet$ & $ \bullet $ & $\bullet$ & $\bullet$ & $\bullet$ & $ \bullet $ & $\bullet$ & $\ldots$\\
         & 0 & & $d_1$ & & $d_1$ & & $d_1$ & & $d_1$ & & $d_1$ & 
    \end{tabular}
\end{center}

\begin{center}
    \begin{tabular}{c c c c c c c c c c c c c}
        $\ldots$ & \hphantom{$\bullet$} & $\bullet$ &  $\hphantom{\bullet}$ & $\bullet$ & $\hphantom{\bullet}$ & $ \bullet $ & $\hphantom{\bullet}$ & $\bullet$ & $\hphantom{\bullet}$ & $ \bullet $ & $\hphantom{\bullet}$ & $\ldots$\\
         & \hphantom{0} & & $\hphantom{d_1}$ & & $\hphantom{d_1}$ & & $\hphantom{d_1}$ & & $\hphantom{d_1}$ & & $\hphantom{d_1}$ & 
    \end{tabular}
\end{center}



\noindent In $S$, the imbedding of multisector $d_2 = 3$ appears as:

\begin{center}
    \begin{tabular}{c c c c c c c c c c c c c}
        $\ldots$ & $\bullet$ & $\bullet$ & $\bullet$ & $\bullet$ & $\bullet$ & $ \bullet $ & $\bullet$ & $\bullet$ & $\bullet$ & $ \bullet $ & $\bullet$ & $\ldots$\\
        & $0$ & & & $d_2$ & & & $d_2$ & & & $d_2$ & & 
    \end{tabular}
\end{center}

\noindent And in $S'$ (odd numbers), the imbedding of multisector $d_2 = 3$ appears as:

\begin{center}
    \begin{tabular}{c c c c c c c c c c c c c}
        $\ldots$ & $\bigodot$ & $\bullet$ &  $\hphantom{\bullet}$ & $\bullet$ & $\hphantom{\bullet}$ & $ \bullet $ & $\hphantom{\bullet}$ & $\bullet$ & $\hphantom{\bullet}$ & $ \bullet $ & $\hphantom{\bullet}$ & $\ldots$\\
         & $(0)$ & & $\hphantom{d_1}$ & $d_2$ & $\hphantom{d_1}$ & & $\hphantom{d_1}$ & & $\hphantom{d_1}$ & $d_2$ & $\hphantom{d_1}$ & 
    \end{tabular}
\end{center}

\noindent In fact, $0$ is not a member of $S'$ (it is not an odd number) so it is represented by $\bigodot$ and labeled $(0)$ to aid the reader. This point must be erased as it does not belong in $S'$. \linebreak When a multisector $d_j$ is imbedded in $S$, the point labeled $0$ must be relabeled $d_j$. \linebreak If these rules are to be followed, then the patterns created by multisector $d_2$ in $S$ and $S'$ appear to be ``identical", i.e., they are both lines of regularly-spaced indistinguishable marks or points stretching endlessly in either direction that look like ... $\bullet$ $\bullet$ $d_2$ $\bullet$ $\bullet$ $d_2$ $\bullet$ $\bullet$ $d_2$ $\bullet$ $\bullet$ ... . Similarly, the imbedding of multisector $d_3 = 5$ in $S$ and $S'$ identically look like ... $\bullet$ $\bullet$ $\bullet$ $\bullet$ $d_3$ $\bullet$ $\bullet$ $\bullet$ $\bullet$ $d_3$ $\bullet$ $\bullet$ $\bullet$ $\bullet$ $d_3$ $\bullet$ $\bullet$ $\bullet$ $\bullet$ ... . Verify this! In general, the imbedding of any multisector in $S$ is identical to its imbedding in $S'$. Hence, $S$ and $S'$ are indistinguishable from one another and $d_1$ is the only multisector that makes this identity. In $S$, the sequence that starts with $0$ can be regarded as starting from some other number such as $d_j$, so $0$ as origin repeats at the primorial $p_1 \times p_2 \times ... \times p_j$, where $p_j = d_j$ is a prime, and these are what Spencer-Brown calls \textit{paradigmatic places} or pp. What distinguishes system $S$ from arithmetic $A$ is that multisector imbeddings in $S$ can be frozen and these patterns can be shifted like a ruler along the number line and put in a \textit{nonparadigmatic place} or np. Spencer-Brown puts these in square segments.





\fbox{\begin{minipage}{32em}
\noindent Every possible arrangement of multisectors in $S$ is matched by a similar arrangement of divisors somewhere in $A$, and the reason for separating the concepts of $S$ and $A$ is to be able to experiment with the multisectors in $S$, by shifting them from place to place, which we cannot do with the divisors in $A$. In other words, I have made in $S$ a working model of $A$, which we can take to pieces and play with the pieces, with the aim of discovering sharp limits to what is possible in $A$ by seeing what is possible, and what is not, in $S$.
\end{minipage}}

\noindent Focusing on the ``imbedding" operation of multisectors, Spencer-Brown uses the word ``strike" with different distinctions of use intended, as a noun, collective noun and operator, collapsed into one by the second canon, \textit{contraction of reference}\footnote{``Let injunctions be contracted to any degree in which they can still be followed" (Chapter 3, \textit{Laws of Form}, Bohmeier Verlag, 2015)}. He writes:

\fbox{\begin{minipage}{32em}
\noindent We can employ the noun \textit{strike} to indicate the appearance of a section of $S$ after the operation of a set of $d$'s. It will appear similar to a section of numbers in $A$ in which a corresponding set of prime divisors have imbedded to make composites, leaving the untouched numbers therein as new primes. (Compare the sieve of Eratosthenes.)
\end{minipage}}


\noindent ``Strike" as a \textit{verb}, or operator, first, is used in reference to the Sieve of Eratosthenes, an ancient algorithm used to generate prime numbers. The following is an illustration of the sieve. Begin with a list of 20 numbers. Omit the number $1$. Begin with the number $2$. Strike all mutiples of $2$ (except itself) off the list. Call the remaining numbers ``holes". After $2$, the first hole is $3$. Use $3$ to strike multiple of $3$ off the list. After $3$, the first hole is $5$. Use $5$ to strike multiple of $5$ off the list, ... The holes that remain are all the primes up to $20$. These determine all primes (holes) up to $20^2$, ...

\begin{center}
    \begin{tabular}{c c c c c c c c c c c c c c c c c c c}
        $2$ & $3$ & $4$ & $5$ & $6$ & $7$ & $8$ & $9$ & $10$ & $11$ & $12$ & $13$ & $14$ & $15$ & $16$ & $17$ & $18$ & $19$ & $20$\\
        \textcircled{$2$} & \textcircled{$3$} & \cancel{$4$} & \textcircled{$5$} & \cancel{$6$} & \textcircled{$7$} & \cancel{$8$} & \cancel{$9$} & \cancel{$10$} & \textcircled{$11$} & \cancel{$12$} & \textcircled{$13$} & \cancel{$14$} & \cancel{$15$} & \cancel{$16$} & \textcircled{$17$} & \cancel{$18$} & \textcircled{$19$} & \cancel{$20$}\\
    \end{tabular}
\end{center}



\noindent In an unpublished manuscript ``Commentary on Appendix 9", Spencer-Brown produces a beautiful method of sieving prime numbers, where $\times$ marks the spot.

\begin{figure}[h!]
            \centering
            \includegraphics[scale=0.55]{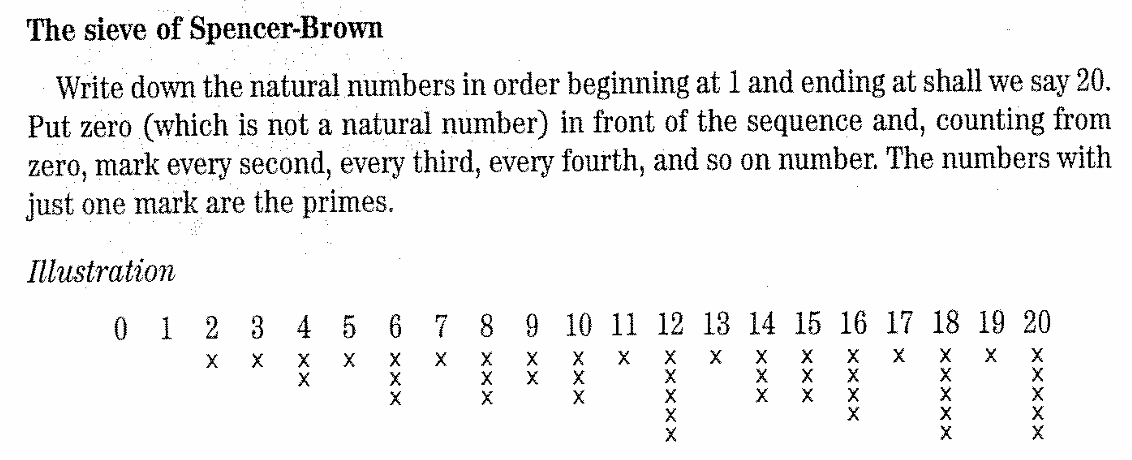}
        \end{figure}

\noindent ``Strike" as a \textit{noun} refers to a representative piece of $S$ or $S'$ when a multisector (or a collection of multisectors) is (are) imbedded in $S$ or $S'$. This representative piece of $S$ in its full form is called the ``full reflexive strike". In the full reflexive strike, we set a 0 point and mark both sides of 0. This yields a four-fold pattern consisting of four quadrants with two ``mirrors" or lines of symmetry located at the primorial (=$p_1\times p_2\times.. \times p_j)$ and half the primorial (which defines the 1st quadrant, or quantum, or half-strike). The following two excerpts are from page 190 of Appendix 8:

\fbox{\begin{minipage}{32em}
\noindent We examine what happens when a set of relatively prime multisectors strike over a consecutive set of points in $S$ or $S'$.
\end{minipage}}



\begin{figure}[h!]
            \centering
            \includegraphics[scale=0.5]{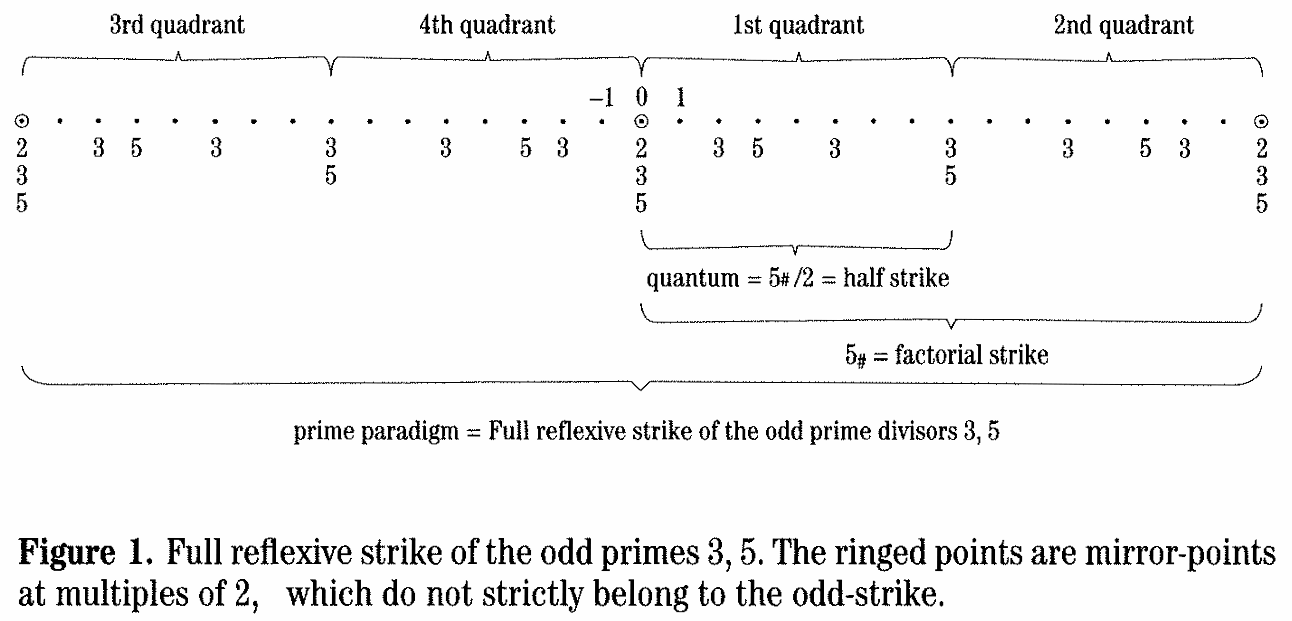}
        \end{figure}

  

\fbox{\begin{minipage}{32em}
\noindent Observe in Figure 1 that we can place inwards-facing mirrors at each end of any quadrant to get the full picture to infinity, so: the quadrant is the \textit{quantum}, or least stretch of space that contains all the visual information we require\footnote{The mirror placings are of particular interest since they must be through the exact location of a point. Thus if these points are euclidean, having position but no size, they will be entirely obliterated by the mirror. To reconstruct them we must allow them an infinitesimal size, so that half the point will appear in the quadrant, and the other half as its image in the mirror.}. Figure 1 shows the repetitive arrangements, and the names I have given them, of a strike of two odd prime decimators or divisors. It can of course be extended to any number of such decimators or divisors.
\end{minipage}}

\noindent ``Strike" as a \textit{collective noun} also refers to the prime factorial, primorial, hash, or prime generator, pg. A definition is provided by Spencer-Brown on page 191 of Appendix 8:

\fbox{\begin{minipage}{32em}
Write $x\#$ (\textit{or $x$!`}), say $x$-hash (\textit{or }$x$\textit{-strike}), for the prime factorial of a nonnegative real number $x$. \hphantom{.} $x\#$ is defined as the factorial of the integer part of $x$, divided by the composite numbers $\leq x$. For example, $\sqrt{80}\#$ $= 8.944...\#$ $= (8!/8.6.4) = 7\#$ $= 210$. Thus the primes $\leq 7$ form the pg of $80$, their strike covers $210$ integers, and their quantum $105$.
\end{minipage}}
%

\thispagestyle{empty}
\section{Sects}

\noindent Spencer-Brown introduces new terminology for the composite numbers that are struck. He calls a sequence of consecutive composite numbers a ``sect". On page 192, he writes:

\fbox{\begin{minipage}{32em}
We shall need one further special term. Mnemonic: imagine a \textit{set} as a Show of Elements Together. Then construct the terms \textit{sect} to denote a Show of Elements Consecutively Together. By a happy chance it extends the ordinary meaning of the word `sect' and perfectly fits its etymology, which is from Middle English \textit{secte}, from Old French, from Latin \textit{secta}, from \textit{sectus}, archaic past participle of \textit{sequi} to follow.
\end{minipage}}

\newpage
\subsection{Odd sects} \hphantom{.}

\noindent From here on out, sects will refer to \textit{odd} sects, sequences of consecutive \textit{odd} composites. We return to our previous example of seg$(5)$. Recall that a square segment, denoted seg$(n)$ is the stretch of dots between squares, and $5$ and $3$ in the third line indicates our convention of writing the first prime that divides the odd number above it.

\begin{center} \begin{tabular}{c c c : c c c c}
        \LARGE $5^2$ & $27$ & $29$ & $31$ & $33$ & $35$ & \LARGE $6^2$\\
        $\square$ & $\bullet$ & $\bullet$ & $\bullet$ & $\bullet$ & $\bullet$ & $\square$\\
        $5$ & $3$ & & & $3$ & $5$ &
    \end{tabular}
\end{center}




\noindent The sect \begin{tabular}{c c}
    $\bullet$ & $\bullet$\\
    $3$ & $5$
\end{tabular} is a sequence of consecutive odd composites in seg$(5)$.\\ 
Another sect we see above is \begin{tabular}{c c}
    $\bullet$ & $\bullet$\\
    $5$ & $3$
\end{tabular}, although this sect is not entirely in seg$(5)$.\\ 
The sects \begin{tabular}{c c}
    $\bullet$ & $\bullet$\\
    $3$ & $5$
\end{tabular} and \begin{tabular}{c c}
    $\bullet$ & $\bullet$\\
    $5$ & $3$
\end{tabular} are maximum odd sects using the pg of $5^2 = \{(2), 3, 5\}$.

\noindent The following are illustrations of the distribution of odd sects in square segments.

\begin{figure}[H]
            \centering
            \includegraphics[scale=0.66]{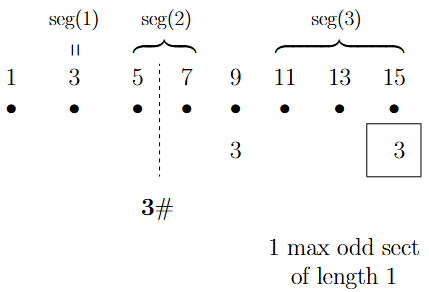}
        \end{figure}

\begin{figure}[H]
            \centering
            \includegraphics[scale=0.65]{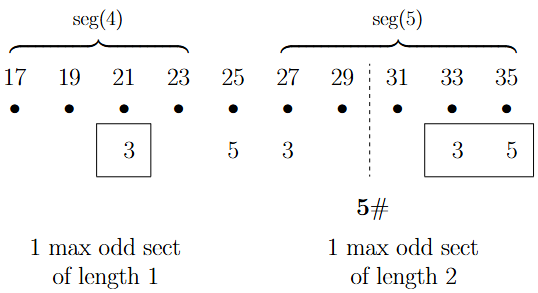}
        \end{figure}

\begin{figure}[H]
            \centering
            \includegraphics[scale=0.66]{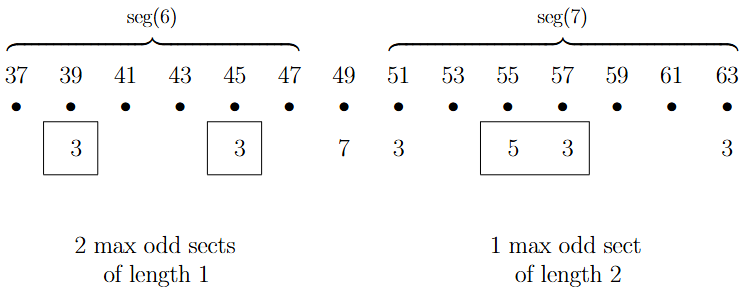}
        \end{figure}

\begin{figure}[H]
            \centering
            \includegraphics[scale=0.64]{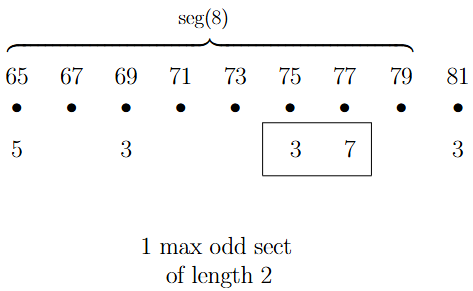}
        \end{figure}

\begin{figure}[H]
            \centering
            \includegraphics[scale=0.64]{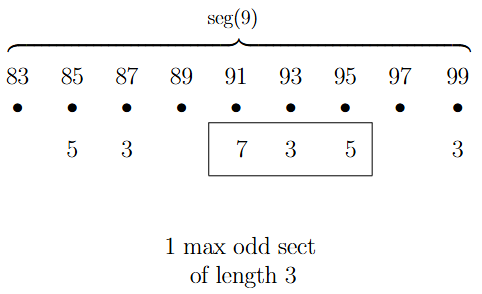}
        \end{figure}

\begin{figure}[H]
            \centering
            \includegraphics[scale=0.66]{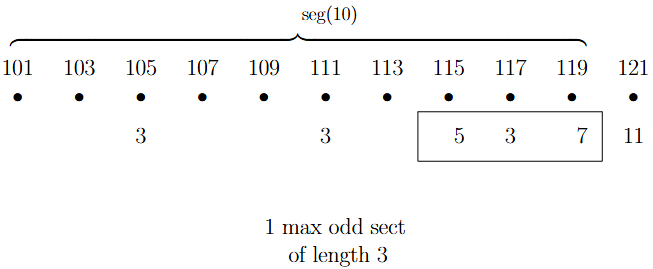}
        \end{figure}

\begin{figure}[H]
            \centering
            \includegraphics[scale=0.68]{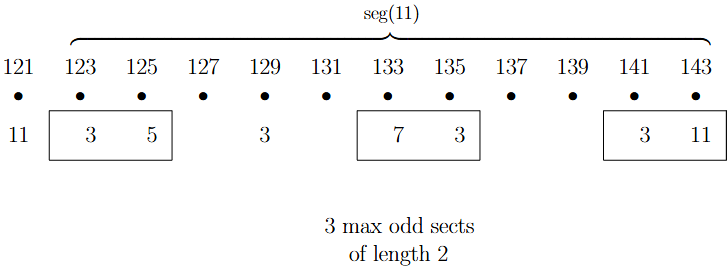}
        \end{figure}

\begin{figure}[H]
            \centering
            \includegraphics[scale=0.68]{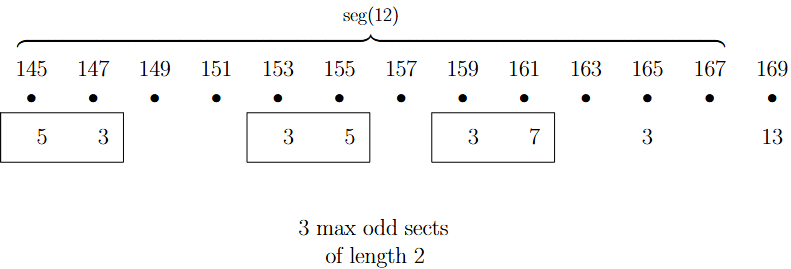}
        \end{figure}

\begin{figure}[H]
            \centering
            \includegraphics[scale=0.68]{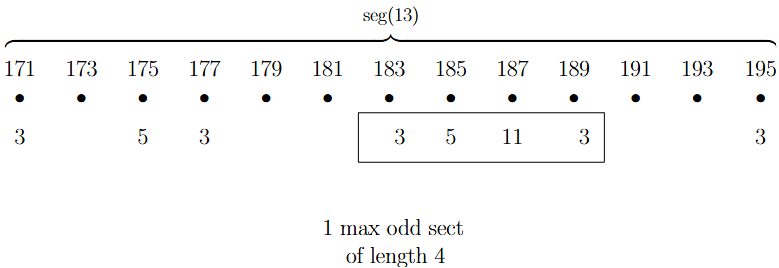}
        \end{figure}

\begin{figure}[H]
            \centering
            \includegraphics[scale=0.62]{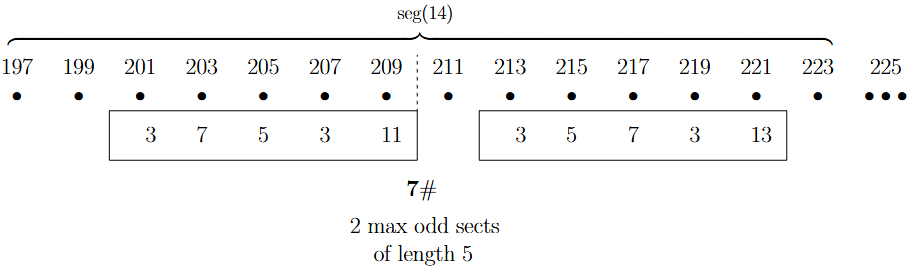}
        \end{figure}









\subsection{What is the maximum sect that can appear in a square segment?} \hphantom{.}

\noindent Spencer-Brown proposes a \textbf{central} question that relates the lengths of maximum sects in square segments to Legendre's conjecture:

\fbox{\begin{minipage}{32em}
What is the maximum set of successive composites that can be generated by the prime divisors in the pg of $n^2$? If the answer comes to less than $2n$, or less than $n$ if we consider the odd numbers only, then we are done, since it will mean that not all the numbers in seg$(n)$ can be struck by these prime divisors, so at least one of the numbers in seg$(n)$ must be prime.
\end{minipage}}

\noindent The answer to Spencer-Brown's question lies in the \textit{prime paradigm}, i.e., Figure 1:


\fbox{\begin{minipage}{32em}
\noindent From Figure 1 we see there are just two basic ways in which a maximum sect of prime divisors $\leq p$ might be constructed.

\noindent 1. By marrying the two reflexive halves of the paradigm of prime divisors as they appear on either side of zero, and thus on either side of any multiple of $p\#$, and occupying the two holes corresponding to $1$ and $-1$ with proper divisors, or (if possible) by some improvement on this recipe, and
    
\noindent 2. By noting the natural sect that begins at $2$ (or using the odd numbers, from $3$ onwards) in the first quadrant, and ends just before the least prime $q > p$.

\noindent These two recipes are together exhaustive: indeed if we are allowed to search the entire quantum of the first quadrant, instead of just the beginning of it, to which we are confined by recipe 2, we shall be certain to find the maximum sect we are looking for, or its mirror image.

\end{minipage}}

%



\noindent Call the sect determined by recipe 1 a \textit{reflexive} or R-sect, and that determined by recipe 2 a \textit{quadrantic} or Q-sect. As we shall see, these recipes can be applied at the origin as well as beyond the origin, namely at multiples of the primorial, $p\#$.

\noindent Spencer-Brown provides the following notation, conflating primes $p$ and prime multisectors $d$, using $p$ and $d$ interchangeably:

\fbox{\begin{minipage}{32em}
Write $p_{\sim j}$ (and similarly $d_{\sim j}$) to indicate the set of all primes (and similarly all prime multisectors) $\leq p_j (\text{or }d_j)$. In $S'$ and the odd number arithmetic $A'$, $p_1$ (or $d_1$) has no place, so $d_{\sim j}$ in $S'$ effectively means $d_2$ through $d_j$.
\end{minipage}}


\newpage

\noindent Let the maximum possible odd R-sect generated by primes multisectors $d_{\sim j}$ be denoted by $\text{m}R(d_{\sim j})$, and the maximum possible odd Q-sect generated by prime multisectors $d_{\sim j}$ be denoted by $\text{m}Q(d_{\sim j})$. Here, Spencer-Brown reuses the notation $\text{m}R(d_{\sim j})$ and $\text{m}Q(d_{\sim j})$ to denote the \textit{lengths} of these maximum possible R-sects and Q-sects, respectively. Hence, Spencer-Brown incorporates both recipes in Lemma 3 as follows:

\fbox{\begin{minipage}{32em}
\noindent \textbf{Lemma 3}:\\
 The maximum possible odd R-sect $\text{m}R(d_{\sim j})$ of consecutive elements in $S'$ that can be struck using all the odd prime multisectors $d_2$ through $d_j$ is\\
 (1) \hphantom{...} $\text{m}R(d_{\sim j}) = d_{j-1} - 1$,\\
 and this sect, counting mirror images as identical, is unique. ...

 \noindent Alternatively, the maximum odd Q-sect will take us to the brink of the next greater prime than $p_j$, yielding the formula\\
  (2) \hphantom{...} $\text{m}Q(d_{\sim j}) = \frac{d_{j+1} - 1}{2}-1$.\\
\end{minipage}}
 

\noindent The observation that $\text{m}R(d_{\sim j}) = d_{j-1} - 1$ was first made by Legendre in \textit{Essai sur la Th\'{e}orie des Nombres}, Quatri\`{e}me Partie, Section IX: "D\'{e}monstration de divers th\'{e}or\`{e}mes sur les progressions arithm\'{e}tiques," pp. 399-406. In addition, Spencer-Brown remarks: ``In fact Lemma 3 identifies the maximum possible sect for the pg of every $n^2$ except in case $n=23$, for which (uniquely) there are six maximum odd sects of $19$". See Appendix 3 for the six unique maximum odd sects generated by prime multisectors $d_{\sim 9}$. A stronger version of Lemma 3 states ``no odd sect can exceed the maximum possible odd Q-sect or R-sect by more than 2".

\noindent Spencer-Brown mentions that proving this stronger version of Lemma 3 is unnecessary. He claims that upon proving Lemma 4, it will become evident that, for $n>5$, the maximum sect that can appear \textit{in the segment} must be increasingly shorter than the maximum sect that can appear \textit{in the quantum}. This assertion is under investigation. The following page illustrates Spencer-Brown's proof of Lemma 3. See section 5.5 of this paper, where we review this algorithm from ground up.

\newpage


\fbox{\begin{minipage}{32em}
\noindent \textbf{Proof of Lemma 3}:

\noindent Suppose we have found Lemma 3 to be true (as it is) for all odd decimators $d_{\sim j}$ when $j = 4$, and for all lesser values of $j$. We have

\begin{spacing}{1.2}
\begin{center}
    \begin{tabular}{c c c c}
        $\bullet$ & $\bullet$ & $\bullet$ & $\bullet$ \\
        $3$ & $5$ & $7$ & $3$
    \end{tabular}
\end{center}
\end{spacing}

\noindent This is the last case where the R-sect and the Q-sect are identical. As a Q-sect, it is merely the strike of odd prime divisors on the positive side of zero, notably over the numbers $3$ through $9$. As an R-sect, we imagine both $3$'s in their paradigmatic places on either side of zero, and arbitrarily choose $5, 7$ to fill the two holes in the strike left by the improper divisors $1$ and $-1$.

\noindent Now consider an induction on $j$. To introduce the next decimator, in this case, $d_5 = 11$, we take $d_3 = 5$ from its nonparadigmatic place between the $3$'s and replace it paradigmatically at both ends of the strike, using $d_5$ to fill the hole in the middle left by $d_3$. In each case we complete the strike, if we can, with the smaller $d$'s already in use. Thus

\begin{center}
    \begin{tabular}{c c c c}
        $\bullet$ & $\bullet$ & $\bullet$ & $\bullet$ \\
        $3$ & $5$ & $7$ & $3$
    \end{tabular}
\end{center}

\noindent \hphantom{...} is followed by

\begin{center}
    \begin{tabular}{c c c c c c}
        $\bullet$ & $\bullet$ & $\bullet$ & $\bullet$ & $\bullet$ & $\bullet$ \\
        $5$ & $3$ & $11$ & $7$ & $3$ & $5$
    \end{tabular}
\end{center}

\noindent \hphantom{...} is followed by

\begin{center}
    \begin{tabular}{c c c c c c c c c c}
        $\bullet$ & $\bullet$ & $\bullet$ & $\bullet$ & $\bullet$ & $\bullet$ & $\bullet$ & $\bullet$ & $\bullet$ & $\bullet$ \\
        $3$ & $7$ & $5$ & $3$ & $11$ & $13$ & $3$ & $5$ & $7$ & $3$
    \end{tabular}
\end{center}

\noindent \hphantom{...} and so on.\\

\noindent We say a multisector is in a \textit{paradigmatic place}, or pp, if it imbeds in a point that can be identified with the associated prime in $A'$. Otherwise we say it is in a nonparadigmatic place np. In any R-sect, the two central multisectors begin by being np. Then, by introducing the next larger multisector, the smallest of them goes from np to pp.\\

\noindent Because at every addition of $1$ to $j$ there is no alternative to what we may do, this validates the induction and proves the uniqueness of the result. QED
 
\end{minipage}}

\fbox{\begin{minipage}{32em}
\noindent \textbf{Lemma 4}:\\
(4) \hphantom{...} $d_{j-1}-1 \geq \frac{d_{j+1}-1}{2}-1$.
\end{minipage}}

\noindent In other words, $\text{m}R(d_{\sim j}) \geq  \text{m}Q(d_{\sim j})$.

\noindent Note that Spencer-Brown does not give a standard arithmetical proof of Lemma 4 but gives a proof of Lemma 4 by \textit{elementals}, which we outline in section 5 of this paper.

\noindent The stronger version of Lemma 4 is $\text{m}R(d_{\sim j}) > \text{m}Q(d_{\sim j})$, for $j>4$, that is, the maximum R-sect is \textit{strictly greater than} the maximum Q-sect for $j>4$.

\noindent Rearranging the terms in Lemma 4 yields a powerful relationship between the ``previous" prime $p_{j-1}$ and the ``successive" prime $p_{j+1}$. Note, $p_j = d_j$.
$$p_{j+1}\leq 2p_{j-1}+1.$$

\noindent An illustration of the above statement for small primes is as follows:\\
$5 \leq 2(2)+1.$\\
$7 \leq 2(3)+1.$\\
$11 \leq 2(5)+1.$\\
$13 < 2(7)+1$. Note: From this point onward, the relationship is a strict inequality.\\
$17 < 2(11)+1$,\\
...

\noindent The reader is encouraged to verify $p_{j+1}\leq 2p_{j-1}+1$ for larger primes!


\noindent What is profound is that Lemma 4 controls how much bigger the $(j+1)$\textsuperscript{th} prime is in relation to the $(j-1)$\textsuperscript{th} prime. In other words, there is a feedback from the $j$\textsuperscript{th} prime to the $(j-1)$\textsuperscript{th} prime and a feed-forward to the $(j+1)$\textsuperscript{th} prime.


\noindent Following the injunctions provided by Legendre in Part 4, Section 9 of \textit{Essai sur la Th\'{e}orie des Nombres} and Spencer-Brown in his proof of Lemma 3, Appendix 8, \textit{Laws Of Form}, we construct what we call the ``sect pyramid" (or the Spencer-Brown sieve or GSB sieve). When viewing this pyramid, compare the maximum odd R-sect of line L with the maximum odd Q-sect of line L+2 to see Lemma 4 in action.

\newpage
\thispagestyle{plain}

\begin{landscape}
\begin{center}
\begin{figure}[h!]
            \centering
            \includegraphics[scale=0.69]{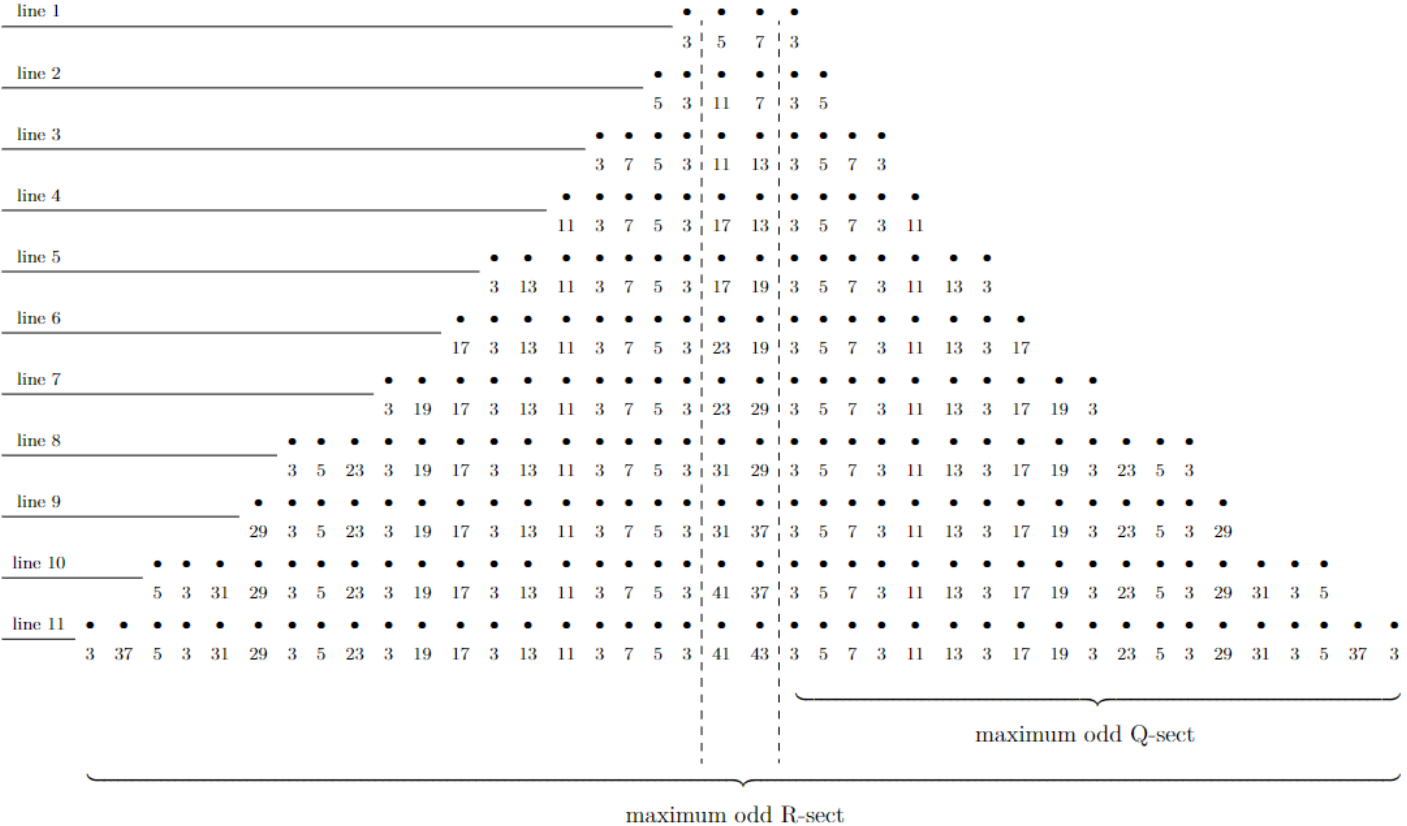}
        \end{figure}
\end{center}
\end{landscape}

\newpage



\subsection{Quadrantic and Reflexive sects} \hphantom{.}

\noindent The maximum odd Q-sect $\text{m}Q(d_{\sim j})$ can be constructed by marking the odd numbers between $1$ and $p_{j+1}$ and replacing each odd composite by its least prime factor.



\noindent Consider $\text{m}Q(d_{\sim 3})$. Marking the odd numbers upto (but not including) $p_4 = 7$, the maximum odd Q-sect $\text{m}Q(p_{\sim 3})$ is \begin{tabular}{c c}
    $\bullet$ & $\bullet$\\
    $3$ & $5$
\end{tabular} (and its mirror image \begin{tabular}{c c}
    $\bullet$ & $\bullet$\\
    $5$ & $3$
\end{tabular}).

\noindent Consider $\text{m}Q(d_{\sim 4})$. Marking the odd numbers upto $p_5 = 11$ (and replacing each odd composite by its least prime factor), the maximum odd Q-sect $\text{m}Q(d_{\sim 4})$ is \begin{tabular}{c c c c}
    $\bullet$ & $\bullet$ & $\bullet$ & $\bullet$\\
    $3$ & $5$ & $7$ & $3$
\end{tabular} (and its mirror image \begin{tabular}{c c c c}
    $\bullet$ & $\bullet$ & $\bullet$ & $\bullet$\\
    $3$ & $7$ & $5$ & $3$
\end{tabular}). 

\noindent In general, it is important to note that the maximum odd Q-sect $\text{m}Q(d_{\sim j})$ appears on the number line mirrored about $p_j$-primorial, i.e. at $p_j\#= p_1 \times p_2 \times \ldots p_j$ and integral multiples of $p_j \#$. This is evident when studying the prime paradigm or full reflexive strike as shown by Spencer-Brown in Figure 1 of Appendix 8. 

\noindent To illustrate the idea that the maximum odd Q-sect $\text{m}Q(d_{\sim j})$ appears on the number line mirrored about $p_j \#$, consider $\text{m}Q(d_{\sim 5})$. Using prime multisectors $d_2$ through $d_5$, the maximum odd Q-sect is $\text{m}Q(p_{\sim 5})$ is \begin{tabular}{c c c c c}
    $\bullet$ & $\bullet$ & $\bullet$ & $\bullet$ & $\bullet$\\
    $3$ & $5$ & $7$ & $3$ & $11$
\end{tabular} (and its mirror image \begin{tabular}{c c c c c}
    $\bullet$ & $\bullet$ & $\bullet$ & $\bullet$ & $\bullet$\\
    $11$ & $3$ & $7$ & $5$ & $3$
\end{tabular}). Note that $p_5$-strike or $p_5\# = 11\# = 2\times 3\times 5\times 7\times 11 = 2310$. On the number line, observe that $\text{m}Q(d_{\sim 5})$ is mirrored about $p_5\#$, which is like an origin:

\begin{center}
\begin{tabular}{| c c c c c c c c}
    \small $1$ & \small $3$ & \small $5$ & \small $7$ & \small $9$ & \small $11$ & \small $13$\\
    $\bullet$ & $\bullet$ & $\bullet$ &  $\bullet$ & $\bullet$ & $\bullet$ & $\bullet$ & $\ldots$\\
    & $3$ & $5$ & $7$ & $3$ & $11$ & 
\end{tabular}
\end{center}

\begin{center}
\addtolength{\tabcolsep}{-2pt}
\begin{tabular}{c c c c c c c c | c c c c c c c c}
    & \small $2297$ & \small $2299$ & \small $2301$ & \small $2303$ & \small $2305$ & \small $2307$ & \small $2309$ & \small $2311$ & \small $2313$ & \small $2315$ & \small $2317$ & \small $2319$ & \small $2321$ & \small $2323$ &\\
    $\ldots$ & $\bullet$ & $\bullet$ & $\bullet$ &  $\bullet$ & $\bullet$ & $\bullet$ & $\bullet$ & $\bullet$ & $\bullet$ & $\bullet$ & $\bullet$ & $\bullet$ & $\bullet$ & $\bullet$ & $\ldots$\\
    & & $11$ & $3$ & $7$ & $5$ & $3$ & & & $3$ & $5$ & $7$ & $3$ & $11$ &
\end{tabular} 
$11\#$
\end{center} 






\noindent On the other hand, the maximum odd $R$-sect $\text{m}R(d_{\sim j})$ is constructed by first placing $p_j$ and $p_{j-1}$ about a ``mirror" (call this a \textit{local} mirror), then appending the maximum odd Q-sect $\text{m}Q(p_{\sim(j-2)})$ and its mirror image on either side of this local mirror. 

\noindent For example, $\text{m}R(d_{\sim 5})$ is constructed by placing $p_5=11$ and $p_4=7$ about a local or internal ``mirror", then appending 
 $\text{m}Q(d_{\sim 3})$ and its mirror image to either side:

\begin{center}
\begin{tabular}{c c c | c c c}
     $\bullet$ & $\bullet$ & $\bullet$ & $\bullet$ & $\bullet$ & $\bullet$\\
    $5$ & $3$ & $7$ & $11$ & $3$ & $5$  
\end{tabular} \hphantom{...} or \hphantom{...} \begin{tabular}{c c c | c c c}
     $\bullet$ & $\bullet$ & $\bullet$ & $\bullet$ & $\bullet$ & $\bullet$\\
    $5$ & $3$ & $11$ & $7$ & $3$ & $5$  
\end{tabular}
\end{center}


\noindent Since the R-sect $\text{m}R(d_{\sim j})$ is composed of two mirrored Q-sects $\text{m}Q(p_{\sim(j-2)})$, the local (or internal) mirror must be located at $p_{j-2}\#$ or some integral multiple of $p_{j-2}\#$. For the example of $\text{m}R(d_{\sim 5})$, the local (or internal) mirror is at $4 \cdot (5\#)$ and $73 \cdot (5 \#)$.

\noindent The maximum odd R-sect $\text{m}R(p_{\sim 5})$ and its mirror image appear on the number line mirrored about $\frac{p_5 \#}{2} = \frac{11\#}{2} = 1155$, which we will call the \textit{global} mirror or half-mirror. This demarcates the first quadrant or the \textit{quantum} (See Figure 1).

\begin{center}
$4\cdot(5\#)$ \hphantom{...................} \hphantom{$\frac{11\#}{2}$} \hphantom{...........................} 
$73\cdot(5\#)$ \hphantom{..}
\addtolength{\tabcolsep}{-1.2pt}
\begin{tabular}{c c c | c c c c : c c c c | c c c}
     \small $115$ & \small $117$ & \small $119$\hphantom{.} & \small \hphantom{.}$121$ & \small $123$ & \small $125$ & \hphantom{......} & \hphantom{......} & \small $2185$ & \small $2187$ & \small $2189$ \hphantom{.} & \small \hphantom{.} $2191$ & \small $2193$ & \small $2195$\\
     $\bullet$ & $\bullet$ & $\bullet$ & $\bullet$ & $\bullet$ & $\bullet$ & $\ldots$ & $\ldots $ & $\bullet$ & $\bullet$ & $\bullet$ & $\bullet$ & $\bullet$ & $\bullet$\\
     $5$ & $3$ & $7$ & $11$ & $3$ & $5$ & & & $5$ & $3$ & $11$ & $7$ & $3$ & $5$
\end{tabular}\\
\hphantom{$4\cdot(5\#)$} \hphantom{..................} $\frac{11\#}{2}$ \hphantom{...........................} 
\hphantom{$73\cdot(5\#)$}
\end{center}




\noindent Next, $\text{m}R(d_{\sim 6})$ is constructed by placing $p_6=13$ and $p_5=11$ about a local or internal ``mirror" and then appending 
 $\text{m}Q(d_{\sim 4})$ on either side. We know that for $\text{m}R(d_{\sim 6})$, the local mirror will be at some integral multiple of $p_4\#$. Using 
 a simple code in Python, we are able to find the first appearance of $\text{m}R(d_{\sim 6})$ centered at $45 \cdot (7\#)$:

\begin{center}
\hphantom{..}$45\cdot(7\#)$\\
\addtolength{\tabcolsep}{-1.2pt}
\begin{tabular}{c c c c c | c c c c c}
$9441$ & $9443$ & $9445$ & $9447$ & $9449$ & $9451$ & $9453$ & $9455$ & $9457$ & $9459$\\
$\bullet$ & $\bullet$ & $\bullet$ & $\bullet$ & $\bullet$ & $\bullet$ & $\bullet$ & $\bullet$ & $\bullet$ & $\bullet$\\
$3$ & $7$ & $5$ & $3$ & $11$ & $13$ & $3$ & $5$ & $7$ & $3$\\
\end{tabular}\\
\end{center}
 
 \noindent The global mirror is at $\frac{p_j \#}{2} = \frac{13\#}{2} = 15015$, so the mirror image of this R-sect is symmetrically placed on the other side of this global mirror, centered at $98 \cdot (7\#)$:

\begin{center}
\hphantom{..}$98\cdot(7\#)$\\
\addtolength{\tabcolsep}{-1.2pt}
\begin{tabular}{c c c c c | c c c c c}
$20571$ & $20573$ & $20575$ & $20577$ & $20579$ & $20581$ & $20583$ & $20585$ & $20587$ & $20589$\\
$\bullet$ & $\bullet$ & $\bullet$ & $\bullet$ & $\bullet$ & $\bullet$ & $\bullet$ & $\bullet$ & $\bullet$ & $\bullet$\\
$3$ & $7$ & $5$ & $3$ & $13$ & $11$ & $3$ & $5$ & $7$ & $3$\\
\end{tabular}\\
\end{center}

\noindent The reader can see evidence of Lemmas 3 and 4 by looking at the illustrations above.

\subsection{A short proof of Theorem 1} \hphantom{.}


\noindent The following proof is slightly different from the proof provided by Spencer-Brown. This proof assumes Lemmas 1A, 2, 3, and 4:

\noindent \textbf{Theorem 1}:\\
There is at least one odd prime between $n^2$ and $(n+1)^2$.


\noindent \textit{Proof}:\\
There are $n$ odd numbers between $n^2$ and $(n+1)^2$ (Lemma 2).

\noindent The least prime factor of any odd composite between $n^2$ and $(n+1)^2$ is $\leq n$ (Lemma 1A). Let $p_1, p_2, p_3, \ldots, p_j$ be all the primes $\leq n$. 

\noindent The longest consecutive stretch of odd composites using $p_2, p_3, \ldots, p_j$ is either $\text{m}R(d_{\sim j})$ of length $p_{j-1} - 1$, or $\text{m}Q(d_{\sim j})$ of length $\frac{p_{j+1} - 1}{2} - 1$ (Lemma 3).

\noindent For $j > 4$, we know that the length of $\text{m}R(d_{\sim j}) > $ the length of $\text{m}Q(d_{\sim j})$ (Lemma 4).

\noindent To prove Theorem 1, the maximum stretch of odd composites between $n^2$ and $(n+1)^2$ cannot be greater than $\text{m}R(d_{\sim j}) = p_{j-1} - 1$. Now $p_{j-1} - 1 < n$, since $p_{j-1} < n$. The consecutive odd composites between $n^2$ and $(n+1)^2$ can only factor into primes up to $p_j$, so the longest length of consecutive odd composites must be $\leq p_{j-1}-1 $, which is $< n$. There are $n$ odd numbers between $n^2$ and $(n+1)^2$ and the longest consecutive stretch of odd composites in this interval is $< n$. Therefore, one of these $n$ odd numbers must be prime.
\qed

\newpage
\thispagestyle{empty}
\section{Generalizing Theorem 1}
\noindent Spencer-Brown's proof of Theorem 1 is subtly different from the proof in section 3.4, but both of these proofs rest on the validity of Lemmas 3 and 4. We will see Spencer-Brown's proof in section 5.1. To prove Lemma 4, Spencer-Brown first generalizes Theorem 1 to Theorem 1A, and proceeds to make powerful generalizations.

\subsection{The first ratchet} \hphantom{.}

\fbox{\begin{minipage}{32em}
\noindent \textbf{Theorem 1A}:\\ 
If $x$ is a real number $\geq 1$, the next prime $> x$ is at a maximum distance of $2\sqrt{x}$ from $x$.
\end{minipage}}

\noindent In other words, Theorem 1A implies that $p_{j+1} \leq p_j + 2\sqrt{p_j}$. This profound statement offers a sharp bound on gaps between primes and it can be sharpened as the primes get larger. This sharpening is called the \textit{ratchet phenomenon}. 

\noindent \textbf{Proposition 1}:\\
Theorem 1A implies Theorem 1.\\
\textit{Proof:} Consider the numbers between $n^2$ and $(n+1)^2$. Let $p_j$ be the largest prime number so that $p_j<n^2$. Let $p_{j+1}$ be the next prime number after $p_j$. Then, by Theorem 1A, $p_{j+1} \leq p_j + 2\sqrt{p_j}$. We have $p_j < n^2$. Therefore $\sqrt{p_j} < n$, whence $p_{j+1} < n^2 + 2n$. Since $(n+1)^2 = n^2 + 2n + 1$, this implies that $p_{j+1} < (n+1)^2$. Hence, there exists a prime $p_{j+1}$ in between $n^2$ and $(n+1)^2$. \qed


\noindent \textbf{Proposition 2 (Andrica equivalent)}:\\ $p_{j+1} \leq p_j + 2\sqrt{p_j}$ implies Andrica's conjecture\footnote{Andrica, D. ``Note on a conjecture in prime number theory". \textit{Studia Univ. Babe\c{s}-Bolyai Math. 1986.} \textbf{31} (4): 44-48.} $\sqrt{p_{j+1}} - \sqrt{p_j} < 1$.\\
\textit{Proof:} Consider, $p_{j+1} \leq p_j + 2\sqrt{p_j}$ implies $p_{j+1} < p_j + 2\sqrt{p_j} + 1$. It follows that $\sqrt{p_{j+1}} < \sqrt{p_j} + 1$. Hence, $\sqrt{p_{j+1}} - \sqrt{p_j} < 1$. \qed

\subsection{Ratchet values and ratchet primes} \hphantom{.}

\noindent Spencer-Brown recasts Theorem 1A in an even more general form:

\fbox{\begin{minipage}{32em}
\noindent $(5) \hphantom{...} p_{j+1} \leq p_j + \frac{1}{k}\sqrt{p_j} \quad \quad (k = \frac{1}{2})$,
\end{minipage}}

\noindent Call $k$ a \textit{ratchet value}. In Theorem 1A, $k=\frac{1}{2}$.

\noindent Now consider, if $p_{j+1} = p_j + \frac{1}{k}\sqrt{p_j}$. This implies $\frac{1}{k}\sqrt{p_j} = p_{j+1}-p_{j}$. Hence, $\frac{1}{k} = \frac{p_{j+1} - p_{j}}{\sqrt{p_{j}}}$ or $k = \frac{\sqrt{p_{j}}}{p_{j+1} - p_{j}}.$ \noindent With this formula in place, one can try values of $k$ for various choices of $j$. Call \textit{ratchet points} (alternatively, \textit{ratchet primes}) those primes $p_{j}$ such that for all $l > j$, $p_{l+1} \leq p_l + \frac{1}{k}\sqrt{p_l}$.

\noindent The first few ratchet primes are provided by Spencer-Brown as $7, 113, 1327, 1669 , …$. The corresponding ratchet values $k = \frac{\sqrt{p_{j}}}{p_{j+1} - p_{j}}$ are
\noindent $k_1 = \frac{\sqrt{7}}{4} = 0.661438...$, then $k_2 = \frac{\sqrt{113}}{14} = 0.7592961...$, then $k_3 = \frac{\sqrt{1327}}{34} = 1.0714120...$.

\fbox{\begin{minipage}{32em}
\noindent We can substitute $k = \frac{\sqrt{7}}{4} = 0.661438...$ when $(5)$ becomes an equality for $i=4$ and remains unequal for all other values of $i$. This depends on knowing that the jump from $7$ to $11$ is the biggest jump in the whole prime paradigm, when compared with the size of the take-off point.

\noindent As soon as this jump has been achieved, $k$ advances to the next ratchet point, which occurs at the jump of $14$ integers from $113$ to $127$ making $k = \frac{\sqrt{113}}{14} = 0.7592961...$ A ratchet point is where $k$ advances to a new value from which it cannot slip back. After $7$, we know that $113$ is the next ratchet point because the jump via $114$ to $126$ is over a maximum R-sect including $11^2$. When this jump has been achieved, $k$ advances to a value $>1$, the next ratchet point being $\frac{\sqrt{1327}}{34} = 1.0714120...$.
\end{minipage}}

\noindent In the handwritten manuscript of \textit{Primes Between Squares}, Spencer-Brown writes:

\fbox{\begin{minipage}{32em}
It is fairly easy to find the ratchet points by empirical observation but to prove they are ratchet points is more difficult. They occur at record jumps between primes, but of course not every record jump indicates a ratchet point. If that were so, the mathematics would be much too easy.
\end{minipage}}

\noindent There is no known formula for the ``next" ratchet point, but entire tables can be made. In the Online Encyclopedia of Integer Sequences (OEIS), there is a sequence of primes associated with Andrica's conjecture produced by Harry J. Smith (OEIS A084974: $7, 113, 1327, 1669, 2477, ...$). It turns out that Smith's prime sequence is apparently identical with the sequence of ratchet primes for our formulation.


\noindent To see how Spencer-Brown's ratchets work in relation to Andrica's conjecture, consider $\sqrt{p_{j+1}} - \sqrt{p_j} < \epsilon$ and examine how large the primes need to be to satisfy this inequality for a given $\epsilon$. If $\sqrt{p_{j+1}} - \sqrt{p_j} < \epsilon$, then $\sqrt{p_{j+1}} < \sqrt{p_j} + \epsilon$. This implies $p_{j+1} < p_j + 2\epsilon \sqrt{p_j} + \epsilon^2$. We can then try ratchet values $\epsilon_j = \sqrt{p_{j+1}} - \sqrt{p_j}$ and look for primes $p_j$ so that for all primes $q > p_j$, the inequality $q' < q + 2\epsilon_j\sqrt{q} + \epsilon_j^2$ holds. Call these special primes $p_j$ \textit{Andrica ratchets}. So far, the Andrica ratchet primes seem to coincide with Spencer-Brown ratchet primes.












\noindent The reader is encouraged to see plots on \textit{Mathematica} that show very strongly how $p_{j+1} \leq p_j + \frac{1}{k}\sqrt{p_j}$ is satisfied for ratchet values $k_r=\frac{\sqrt{p_r}}{p_{r+1}-p_r}$, where $p_r=7, 113, 1327, ...$.



\noindent Another \textit{Mathematica} experiment is to choose a $k \neq k_r$, that is a value of $k$ that does not correspond to a ratchet prime. The plot behaves chaotically, but then beyond a certain point the inequality $p_{j+1} \leq p_j + \frac{1}{k}\sqrt{p_j}$ holds. This is a numerical illustration of Spencer-Brown's Theorem 1B.

\fbox{\begin{minipage}{32em}
\noindent \textbf{Theorem 1B}:\\
The size of $k$ can be as large as we please, and $ p_{j+1} \leq p_j + \frac{1}{k}\sqrt{p_j}$ will be true of all numbers from some $j$ upwards.
\end{minipage}}

\noindent The following is a discrete plot of $p_j + \frac{1}{k}\sqrt{p_j} - p_{j+1}$, for arbitrary $k\approx 2.61008$.\linebreak Note how from $j\approx 3400$ upwards, $p_j + \frac{1}{k}\sqrt{p_j} - p_{j+1} > 0$.

\begin{figure}[h!]
            \centering
            \includegraphics[scale=0.34]{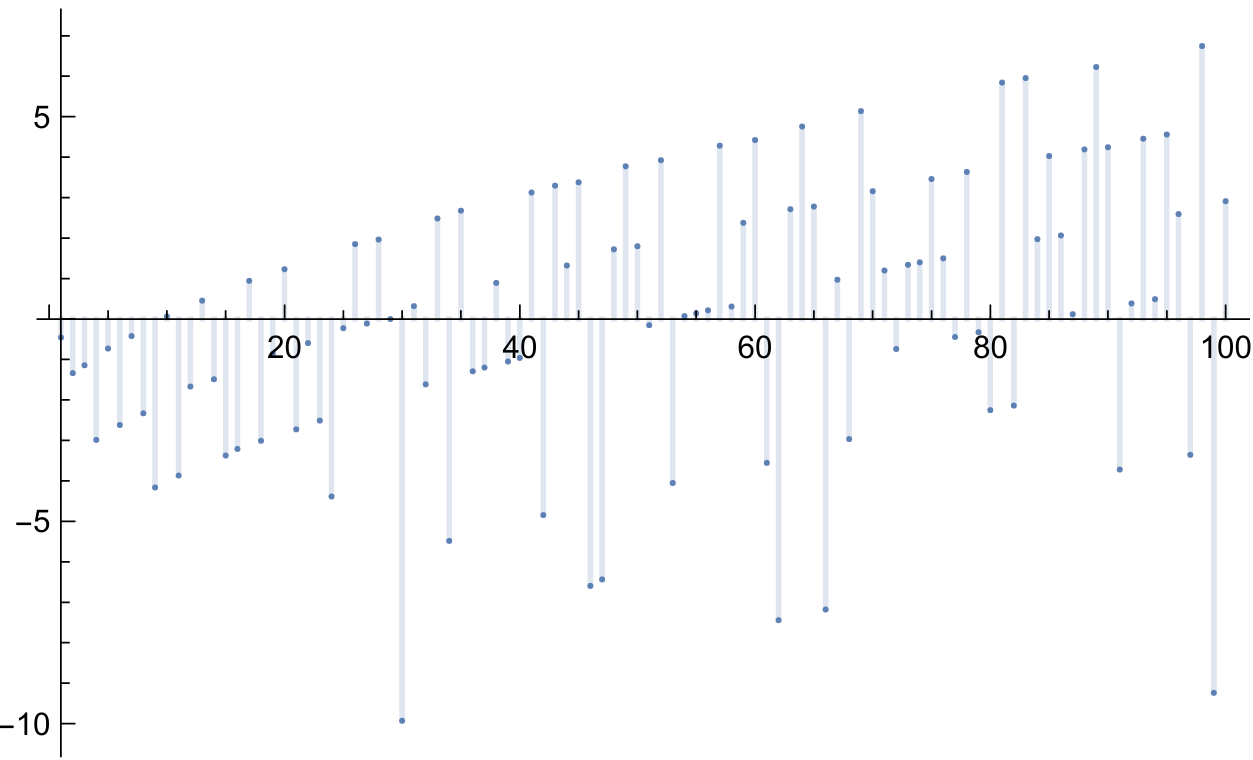}
        \end{figure}

\newpage
\hphantom{.}
\vspace{1em}

        \begin{figure}[h!]
            \centering
            \includegraphics[scale=0.34]{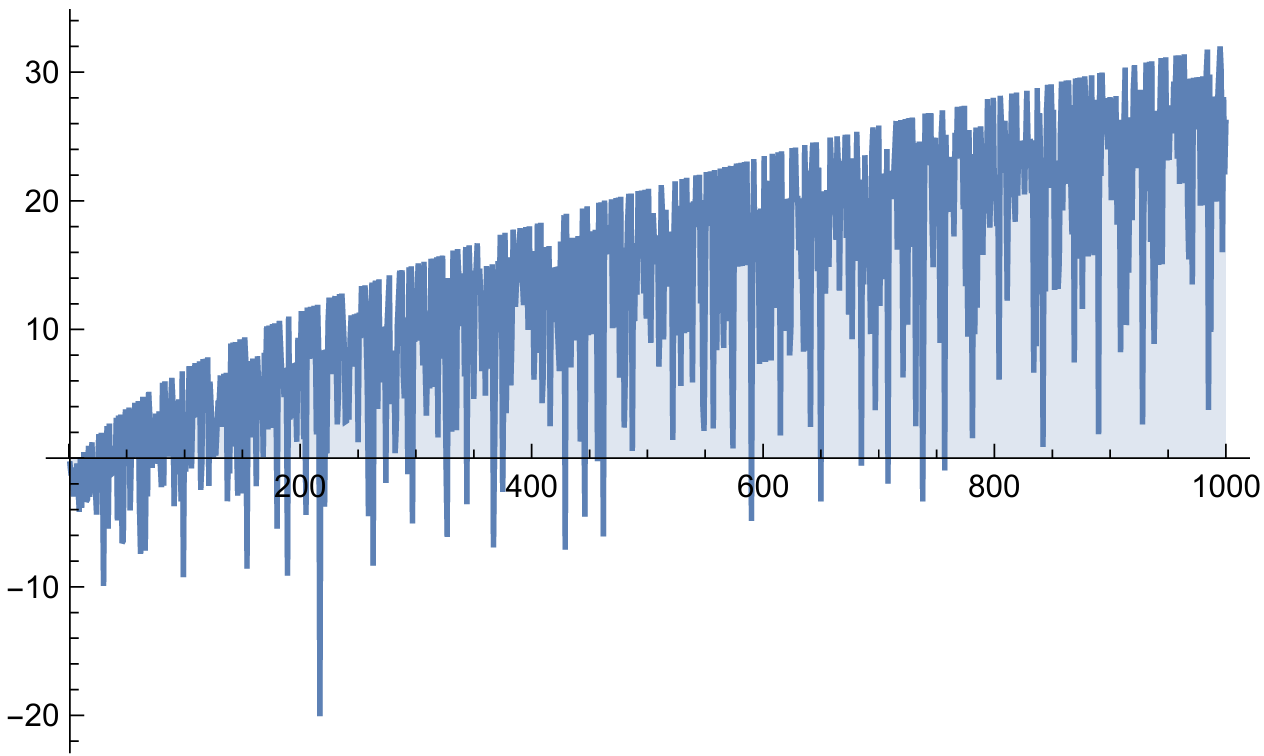}
        \end{figure}

\begin{figure}[h!]
            \centering
            \includegraphics[scale=0.38]{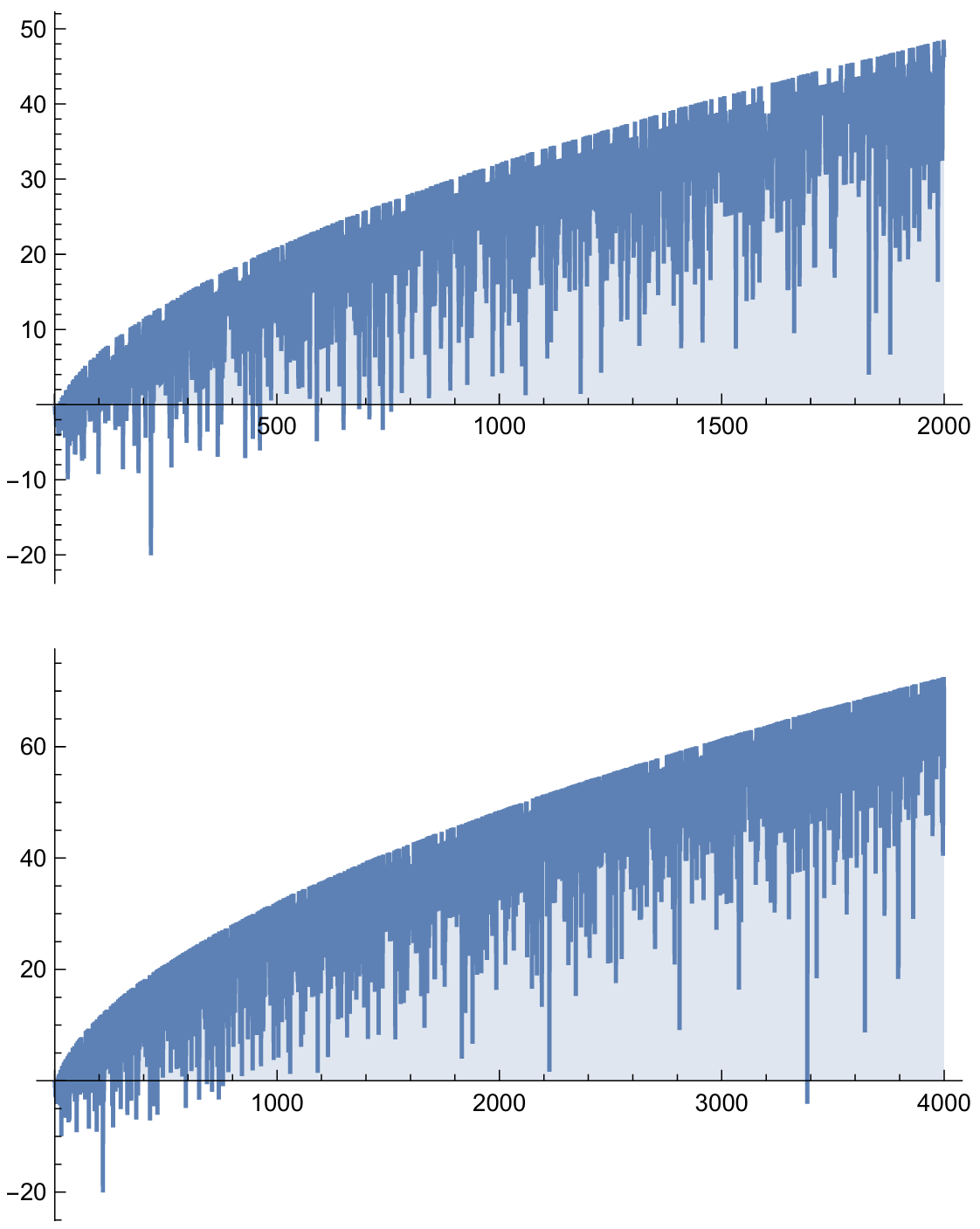}
        \end{figure}

\begin{figure}[h!]
            \centering
            \includegraphics[scale=0.38]{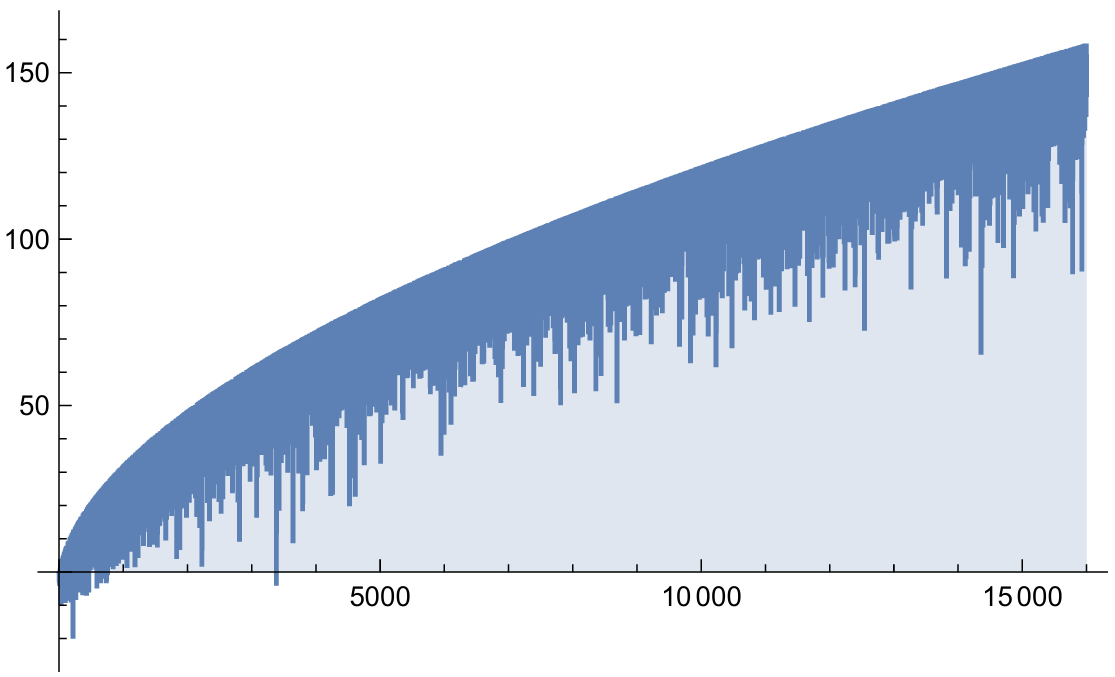}
        \end{figure}


\newpage
\subsection{Haploid and Diploid numbers}\hphantom{.}

\noindent Spencer-Brown then provides increasingly powerful and insightful statements about the distribution of primes between squares. The first of these is a corollary\footnote{See Baker, Harman, Pintz (2000) for a proof of the theorem ``There exists $x_0$ such that for all $x>x_0$, the interval $[x-x^{0.525},x]$ contains a prime numbers". In Spencer-Brown's personal copy of A. E. Ingham's \textit{The Distribution of Prime Numbers}, we find the following handwritten conjecture: ``p 8 \hphantom{.} RH $\equiv$ or implied by $p_{n+1}-p_n < p_n^{1/2}$ for some $n$ onwards".}, which states that the nearest prime neighbors of any given prime $p_j \geq p_{32}$ lie in the intervals $(p_j - \sqrt{p_j},p_j)$ and $(p_j,p_j+\sqrt{p_j})$. 


\fbox{\begin{minipage}{32em}
\noindent \textbf{Corollary}:\\
If $j \geq 32, (p_{j+1} - p_j)$ and $(p_j - p_{j-1})$ are both less than $\sqrt{p_j}$
\end{minipage}}

\noindent Note, the ratchet prime $113$ is the 31\textsuperscript{st} prime, $p_{31}$. This corresponds to the ratchet value of $k = \frac{\sqrt{p_{31}}}{p_{32}-p_{31}} \approx 0.7592961$. The next ratchet prime $1327$ is the 217\textsuperscript{th} prime, $p_{217}$. This corresponds to the ratchet value of $1.0714120$. We note that somewhere between $p_{31}$ and $p_{217}$, the value $k$ crosses a value of $1$ and grows larger (Theorem 1B), implying that $\frac{1}{k}$ crosses a value of $1$ and grows smaller. Hence, $p_{j+1} \leq p_j + \frac{1}{k}\sqrt{p_j}$ becomes $p_{j+1} \leq p_j + \sqrt{p_j}$ from here on out. This provides some insight into the Corollary.


\noindent We also make reference to Spencer-Brown's unpublished ``Haploid and Diploid Numbers", in which he writes ``The largest haploid number is $127$. All numbers greater than $127$ are diploid". A haploid number $n$ is a number that has only one prime neighbor in either $(n-\sqrt{n},n)$ or $(n,n+\sqrt{n})$. A diploid number has a prime neighbor in both $(n-\sqrt{n},n)$ and $(n,n+\sqrt{n})$. Here, Spencer-Brown is making reference to the Corollary, specifically the ratcheting phenomenon that occurs beyond $127$, the 32\textsuperscript{nd} prime, after which $p_{j+1} \leq p_j + \sqrt{p_j}$.


\noindent The following is a discrete plot of $k-1 = \frac{\sqrt{p_j}}{p_{j+1}-p_j} - 1$ on \textit{Mathematica}. Note that the last time $\frac{\sqrt{p_j}}{p_{j+1}-p_j} < 1$ is at $j=31$. This is a numerical illustration of the Corollary.

\begin{figure}[h!]
            \centering
            \includegraphics[scale=0.9]{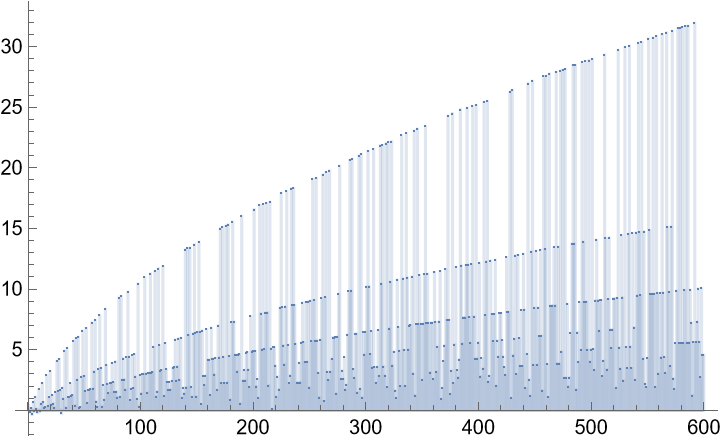}
        \end{figure}


\subsection{Half segments} \hphantom{.}

\noindent Recall that seg$(n)$ is the stretch of numbers between $n^2$ and $(n+1)^2$. A \textit{half segment} is the stretch of numbers between $n^2$ and $n(n+1)$, or between $n(n+1)$ and $(n+1)^2$.

\noindent Note that seg$(1)$ is the stretch of numbers between $1$ and $4$. Halving this segment illustrates the existence of one prime in each half segment, $2$ and $3$. Similarly, note that seg$(2)$ is the stretch of numbers between $4$ and $9$. Halving this segment illustrates the existence of one prime in each half segment, $5$ and $7$. Demonstrate the existence of one prime in each half segment for seg$(3)$ and seg$(4)$. Observe seg$(5)$ below. In this case, the primorial $2\times 3\times 5 = 30$ is the dashed line that splits seg$(5)$ into half segments, illustrating the existence of one prime in each half segment, $29$ and $31$.

\begin{center} \begin{tabular}{c c c : c c c c}
        \LARGE $5^2$ & $27$ & $29$ & $31$ & $33$ & $35$ & \LARGE $6^2$\\
        $\square$ & $\bullet$ & $\bullet$ & $\bullet$ & $\bullet$ & $\bullet$ & $\square$\\
        $5$ & $3$ & & & $3$ & $5$ &
    \end{tabular}
\end{center}

\noindent On page 200 of Appendix 8, Spencer-Brown writes: ``I will close with a more-difficult theorem, that I first proved in the summer of 1998". He finally concludes with the following conjecture:

\fbox{\begin{minipage}{32em}
\noindent \textbf{Theorem 2}:\\
For all natural squares, there is a prime in every half segment\footnote{Oppermann's (hitherto unproved) conjecture of 1862 is incomplete, because its wrong notation disallows the listing of $2$ between $1^2$ and $(1+\frac{1}{2})^2$. See Dickson's \textit{History} Vol 1 p 435. I can thus claim that this delightful theorem is mine.}.
\end{minipage}}






\noindent Ulam's spiral\footnote{Stein, M. L. Ulam, S. M. Wells, M. B. ``A Visual Display of Some Properties of the Distribution of Primes". Amer. Math. Monthly 71, 516-520. 1964.} has a variation called the ``counterdiagonal of squares" where every completion of a square grid of numbers ends with a square number. This variation clearly demarcates successive square segments as a pair of adjacent edges of the square grid. As seen below, each edge of this grid is a half segment. Theorem 2 states that there exists a prime in every edge of this counterdiagonal of squares:



\begin{figure}[h!]
            \centering
            \includegraphics[scale=0.57]{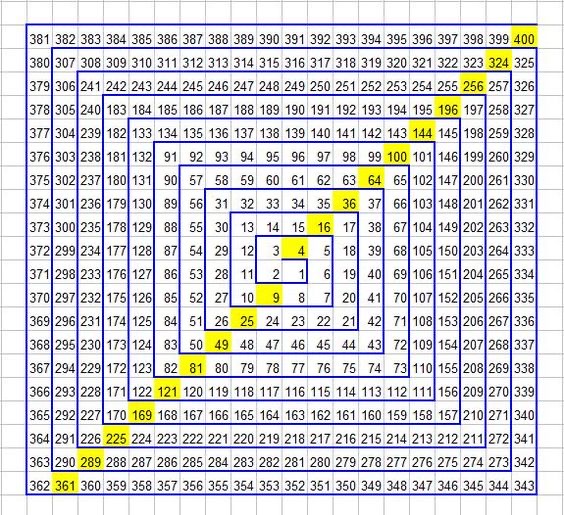}
        \end{figure}


\newpage
\thispagestyle{empty}
\section{Elementals}
\subsection{Spencer-Brown's proof of Theorem 1} \hphantom{.}

\noindent Spencer-Brown's proof of Theorem 1 in Appendix 8 employs a different form of reasoning. A brief outline of his proof is as follows.


\fbox{\begin{minipage}{32em}
\noindent \textit{(Spencer-Brown's proof of Theorem 1)}\\
If Lemma 4 is true, as we are supposing, then Lemma 3 identifies the maximum possible sect that might appear in seg$(n)$, as a result of some strike of all the prime multisectors $\leq n$, and the proof of Theorem 1 can be completed without further difficulty. We must show that the inequality 
$$n > \text{m}R(d_{\sim j})$$
is always true when $j$ is the ordinal number of the greatest prime in the pg of $n^2$.

\noindent When $n=1$ or $n=2$, we cannot use just the odd numbers. There is no prime in the pg of $1^2$, so all the integers in seg$(1)$, $2$ and $3$, must be prime, and $3$ is odd. There is only one prime $2$, in the pg of $2^2$, so all the odd numbers, $5$ and $7$, in seg$(2)$ must be prime. For larger $n$ we can consider the odd numbers only. There can be no odd prime nearer to $n$ on the small side than $n-2$, and since by Lemma 3, $\text{m}R(d_{\sim j}) = d_{j-1}-1$, the maximum odd sect must be of length $p_{j-1}-1$, i.e., is always short of $n$ by at least $3$.\footnote{The special case of $n=23$ makes to difference, since $23$, not being a lesser prime twin, has no effective strike in seg$(n)$.} So the maximum successive strike of odd numbers by the odd primes in the pg of $n^2$ is at most $n-3$, so at least one of the $n$ odd numbers between $n^2$ and $(n+1)^2$ is prime. QED
\end{minipage}}

\noindent Note that the proof of Theorem 1 in section 3.4 is similar to the above proof above. Additionally, Spencer-Brown points out a circular relationship between Theorem 1, Lemma 3, and Lemma 4. This is illustrated as follows and described in detail, shortly.
\begin{center}
Lemma 4 $\implies$ Lemma 3 $\implies$ Theorem 1 $\implies$ Lemma 4 ...  
\end{center}

\newpage

\noindent After his proof of Lemma 3, he writes:

\fbox{\begin{minipage}{32em}
Furthermore, we are now aware that Lemma 3 will always indicate the maximum\footnote{Because for $n>5$ the pg divisors cannot all be pp in the segment. For $n>5$ the maximum sect in the segment will be shorter by an increasing margin than the maximum R-sect, provided the inequality in Lemma 4 is true by an increasing margin, as I shall it must be.} possible sect that might appear in seg$(n)$, provided Lemma 4 is true, in other words, provided the primes in the paradigm are not too far apart.
\end{minipage}}





\noindent According to Spencer-Brown in the above excerpt, Lemma 3 is contingent on Lemma 4. And in subsequent excerpts, he states that Lemma 4 is contingent on Theorem 1. To prove Lemma 4, Spencer-Brown constructs Theorem 1A by ``substituting a continuous variable $x$ for the stochastic variable $n$ in Theorem 1" and introduces fake prime tables, Table 1 and Table 2.



\noindent Starting with the idea of ``not too far apart", Spencer-Brown considers the extreme case that the primes $p_j > 5$ are distributed in such a way that the next prime, $p_{j+1}$, would be of the form $2p_j - 3$. This is evidently not how the primes are distributed in reality but this experiment is used by Spencer-Brown to show that Bertrand's postulate (Tchebychef's theorem\footnote{For $n>3$, there is a prime between $n$ and $2n-2$.}) is ``misleadingly weak" and the resulting calculations for maximum R and Q sects cannot be possible. On page 195 of Appendix 8, he writes:


\fbox{\begin{minipage}{32em}
\noindent Suppose we could show that Theorem 1A is true in all such instances of $x$ after a certain point. How would this affect Lemma 4? First we construct a table of maximum distances between successive primes that would render Lemma 4 barely true.
\end{minipage}}

\noindent Table 1 in Appendix 8 is a fake prime table. Here, Spencer-Brown follows a line of reasoning similar to his earlier assumption of Bertrand's postulate, now assuming that Lemma 4 is ``barely true". Lemma 4 states that $\text{m}R(d_{\sim j}) \geq  \text{m}Q(d_{\sim j})$. By Lemma 4 being ``barely true", Spencer-Brown means $\text{m}R(d_{\sim j}) = \text{m}Q(d_{\sim j})$, i.e., he assumes that $p_{j+1} = 2 p_{j-1} + 1$. Using this relationship, he determines the next ``prime", as shown above, in column 2.


\begin{figure}[h!]
            \centering
            \includegraphics[scale=0.6]{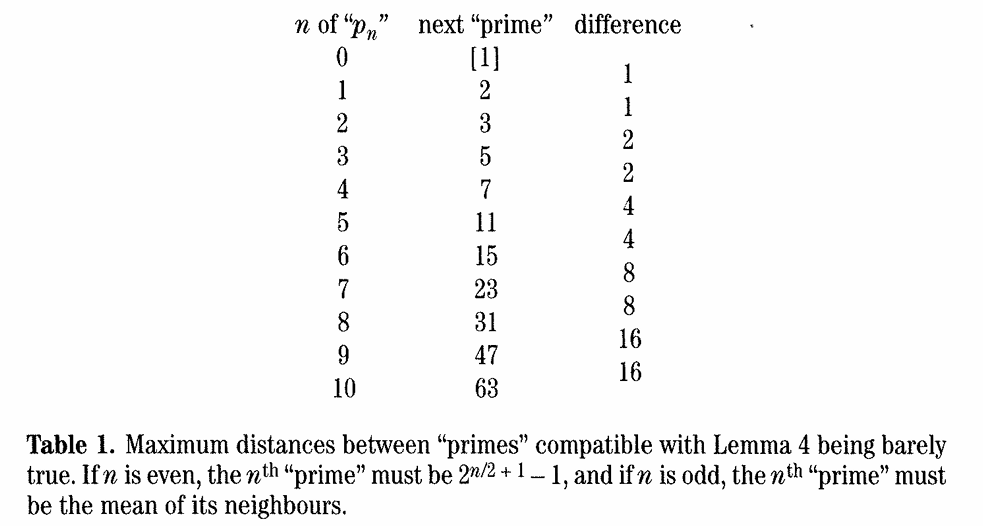}
        \end{figure}


\noindent Table 2 in Appendix 8 is also a fake prime table.\footnote{Note: The caption of Table 2 on page 196 of \textit{Laws Of Form} (Bohmeier Verlag, 2015) should be corrected to $p_{j-1}$, $p_j$, and $p_{j+1}$.} Spencer-Brown defines ``successive" ``prime" here as the nearest odd number within the bare limits of Theorem 1A. In other words, a ``successive" ``prime" is the nearest odd number within two square-roots of the actual prime $p_{j-1}$ in column 1. To read this table, start with $p_{j-1}$, which is the real prime. Take $11$, for example, then the next prime is $11 + 2\sqrt{11} \approx 17$, and the next prime is $17+2\sqrt{17}\approx 25$. 


\begin{figure}[h!]
            \centering
            \includegraphics[scale=0.6]{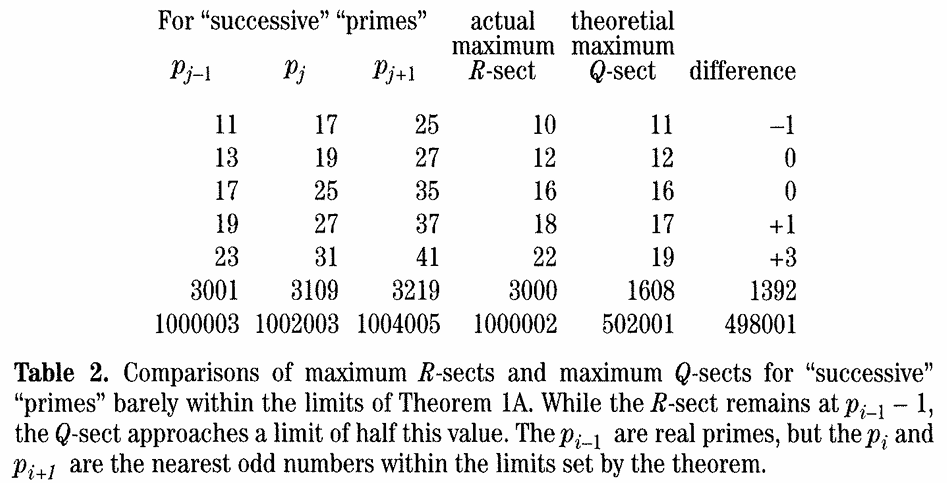}
        \end{figure}






\noindent Using Table 1 and Table 2, Spencer-Brown provides the following analysis on page 197:

\fbox{\begin{minipage}{32em}
\noindent Comparing where the tables become comparable, we see that, for $p_{j-1} = 23$, Table 2 makes a span $23, 31, 41$, while Table 1 makes a wider span $23, 31, 47$, so Table 2 from here on will always give values within the limits of Lemma 4 being true. In fact it begins to happen at $p_{j-1} = 13$, and from $p_{j-1} = 19$ the maximum R-sect is always greater than the maximum Q-sect. Thus if Theorem 1A is true, it locks Lemma 4 into inequality from $n=19$ onwards. This in turn ensures that Lemma 3 must indicate the longest possible sect of the multisectors in the pg of $n^2$ that can appear in seg$(n)$. And if Lemma 3 from then on always indicates the longest possible sect in seg$(n)$ of the multisectors corresponding to the prime divisors in the pg of $n^2$, then Theorem 1A, and with it Theorem 1, must be true from then on. And if Theorem 1 is true from then on, Lemma 4 is true, and therefore irrelevant, from then on, therefore Lemma 3 indicates the maximum sect in seg$(n)$ from then on, and so on \textit{ad infinitum}.\\

\noindent It is a simple matter to ascertain that Theorem 1 is true for $n$ up to $19$, and is true for the intervening values of $n$ to $23$, and so on, so that Lemma 4 cannot be false for these values either. So Theorem 1, by proving its own lemma, has rendered itself true for all $n$, and no other argument is either necessary or relevant.
\end{minipage}}


\noindent ``In short, Theorem 1, like the root in Newton's formula\footnote{Newton's formula or the Newton-Raphson method (named after the method found in Isaac Newton's \textit{Method of Fluxions} and Joseph Raphson's \textit{Analysis Aequationum Universalis}) is an iterative method for finding the roots of an arbitrary function $f(x)$. Starting with an initial guess $x_0$, a better approximation $x_{n+1}$ is obtained by using the iterative formula $x_{n+1} = x_n - \frac{f(x_n)}{f'(x_n)}$, where $f'(x_n)$ is the derivative of $f(x)$ at $x_n$. The process is repeated until the approximations become close to the actual root.}, hunts out its own truth" (page 198, \textit{Laws Of Form}, Bohmeier Verlag, 2015). In a footnote to this statement, Spencer-Brown writes: ``And so do Lemmas 3 and 4, since each elemental of an interlocking set is of equal status. Thus we need no longer consult Table 2 to see experimentally that $\frac{\text{mQ}d_{\sim j}}{\text{mR}d_{\sim j}} \to \frac{1}{2}$ as $j \to \infty$, because the interlocking truth of Lemma 4 with Theorem 1A now assures us \textit{mathematically} that this must be so".

\noindent Thus, Spencer-Brown defines \textit{elementals}, a recursive form of mathematical reasoning whose algebraic representation is a ``valid equation of the second degree, with real, i.e., self-confirming, roots".\footnote{The quadratic equation $x^2=ax+b$ can be reformulated as $x=a+\frac{b}{x}$. This can be further reformulated as a continued fraction: $x = a+\frac{b}{a+\frac{b}{a+\frac{b}{a+...}}}$. This is an instance of a valid equation of second degree with self-confirming roots. If the discriminant of this quadratic equation is non-negative, the roots are real numbers. However, if the discriminant of this quadratic equation is negative, the roots are complex numbers. This continued fraction is the root of the quadratic equation.} The following two excerpts are from pages 197-198.

\fbox{\begin{minipage}{32em}
I coined the new noun \textit{elemental} to denote any of a set of consequential theorems that dovetail with each other in this way, so that they all stand or fall together. To prove that any of them is true, it is necessary first to suppose that one of them is, whereafter the resulting truth of the others confirms the truth of the initial one. It is a kind of contrary to an ordinary indirect proof, where we first have to suppose the initial proposition is false, and get a contradiction. Here we can suppose it is true, and get a confirmation.
\end{minipage}}

\noindent If one walks through the following argument using \textit{elementals} very carefully, letting the structure count, it has the structure of the modulator function E4, as introduced in Chapter 11 of \textit{Laws Of Form}. We will discuss this in detail shortly.



\fbox{\begin{minipage}{32em}
\noindent I suggest the reader pause for a moment to consider the astonishing nature of this argument. We have three theorems, Lemma 3, considered as indicative of the longest sect that can appear in $\text{seg}(n)$, Lemma 4, and Theorem 1. They are all different, i.e., they are not equivalent, and each can vary within certain limits. None of them can be proved unless we assume first that one of them is true. We cannot prove that Lemma 3 produces the maximum sect that can appear in $\text{seg}(n)$ unless Lemma 4 is true, we cannot prove Lemma 4 is true unless Theorem 1 is true, and we cannot prove Theorem 1 is true unless Lemma 3 produces the longest sect that can appear in $\text{seg}(n)$. All these theorems stand or fall together. If one of them is true, they must all be true, and if one of them is false, they may all be false. Yet there is no way of proving any one of them is true without knowing that another of them is also true.
\end{minipage}}

\noindent An illustration of his approach using elementals is provided in the diagram below. Here, each arrow represents implication. It is important to note that this diagram only captures one chain of implications, and the structure of relationships between these statements, when allowed to count, is seen to follow a modulator relationship.

\begin{figure}[h!]
            \centering
            \includegraphics[scale=0.6]{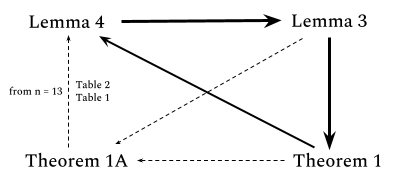}
        \end{figure}





\noindent In the ``sect pyramid", we can see that the maximum odd R-sect in line 1 and the maximum odd Q-sect in line 3, are equal. This corresponds to the maximum odd sects generated by the first $j$ primes, $j=4$. From $j=5$ onward (compare R-sect of line 2 and Q-sect of line 4, and so on, ...), the reflexive sect is greater than the quadrantic sect (Lemma 4) and new holes are introduced. Spencer-Brown uses Table 1 and Table 2 to provide numerical evidence of the first ratchet from $p_4 = 7$ to $p_5 = 11$. The following section, taken from page 200 of Appendix 8 and the handwritten edition of \textit{Primes Between Squares}, points out how beyond this ratchet point, $^\circ \pi$ seg $(n) \geq 3$.

\fbox{\begin{minipage}{32em}
By associating primes with empty spaces in the strike, we can produce calculated theorems much stronger than Theorem 1. Because $n\#$ grows so much faster than $n^2$, we see that $n=5$ is the last time that $n\#$ or any multiple of it can occur within seg$(n)$ with all the pg divisors paradigmatically placed. In this segment we have exactly two primes, $29$ and $31$, corresponding to the two unoccupied central points in the R-strike of $2, 3, 5$, all of which are pp. From here on at least one of the pg divisors will have to be np, making at least one extra hole in the strike. [Therefore when $n$ is a natural number $\geq 6$, $^\circ \pi$ seg $(n) \geq 3$. QED (handwritten)] In practice each np divisor leaves an average of about one extra hole, so an approximate answer for the number of primes in a segment is 2 plus the number of np divisors in the pg.
\end{minipage}}

\begin{figure}[H]
            \includegraphics[scale=0.23]{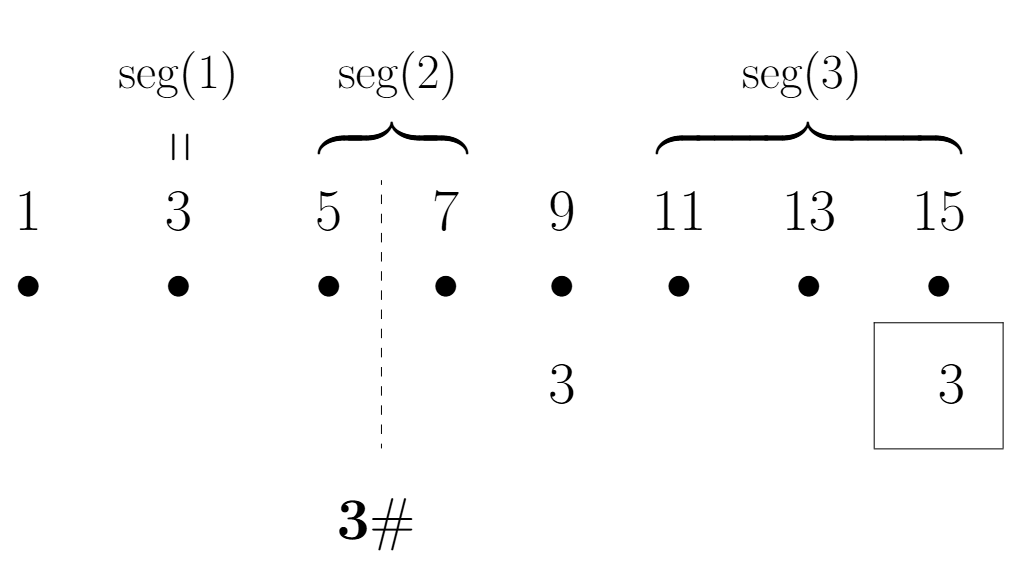}
             \includegraphics[scale=0.23]{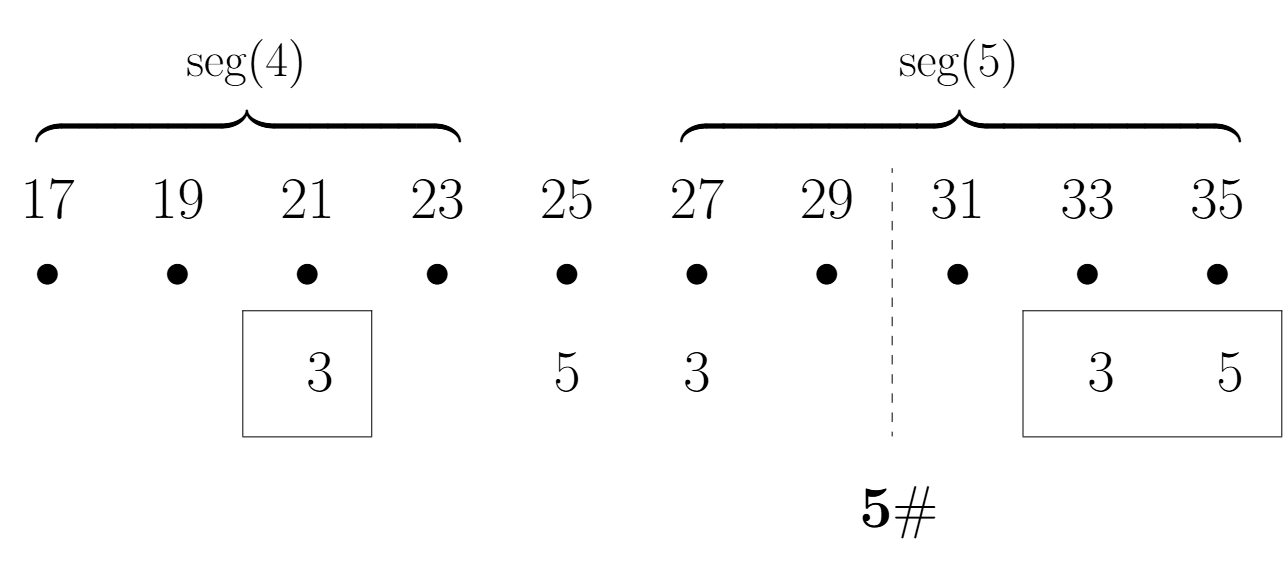}
        \end{figure}

\vspace{-2em}

\noindent Recall that $n\#$ refers to the product of all primes up to $n$. In the figure above, observe that $5\#$ is in seg$(5)$ and this is the last time that $n\#$ can occur within seg$(n)$. Note that $7\#$ is in seg$(14)$, and $n\#$ will always lie outside of seg$(n)$ for all $n>5$. This key excerpt on page 200 of Appendix 8 highlights exactly where the E4 modulator cycle becomes a spiral, starting with the first ratchet. If it were a true modulator, it would repeat itself modulo $4$ (resulting in Table 1). The following section reviews the modulator E4 closely in this context.

\subsection{Modulators}  \hphantom{.}

\noindent Chapter 11, \textit{Laws of Form}, is referenced\footnote{``The reader may refer to Appendix 8, where I use an argument of degree 2 to prove that there are primes between squares" (Original Footnote, page 81).} in both printed and handwritten editions of Appendix 8, \textit{Primes Between Squares}. The ``modulator" or ``reductor" in Chapter 11, called E4, is a re-entering circuit using NOR gates that takes an input 1010 and returns an output 1100 or 0011, i.e., halves the input frequency, over time. In chapter 11, Spencer-Brown points out that imaginary values are important for the structures of modulators. What Spencer-Brown means by an ``imaginary value" is a ``transitional state" that maintains a value just long enough to ensure a desired outcome. For modeling purposes, we assume that there is some time delay across the mark but lines transmit information instantaneously. The modulator has stable states that when perturbed by a single bit input make a transition to a new (stable) state. The modulator is designed in such a way that the transition state is unique. With this determinate structure of transitions, the modulator is an entity that can produce a temporal pattern but is itself outside of time, hence, eternal.\footnote{``The reader's attention may be drawn, for example, to the parallel accounts of the emergence of time, i.e. the statements of what we have to do to construct an element that doesn't exist in any of the five orders of eternity" (James Keys, \textit{Only Two Can Play This Game}, pg. 108, 123-132).}





\noindent There are two ways of understanding modulators -- one is by tracking what is happening locally at each marker when the modulator is in a balanced state\footnote{For a lucid discussion on balanced states and transitions of E4, see ``Modulators and Imaginary Values" by Kauffman (Chapter 7 of \textit{Laws Of Form - A Fiftieth Anniversary}, 2020).} and another is by tracking the waveform that emerges globally at each marker. The following are diagrams of four stable states. Each state indicated after the vertical line undergoes a transition to the next one. We do not here go into the details of the transitions.

\begin{figure}[h!]
            \centering
            \includegraphics[scale=0.17]{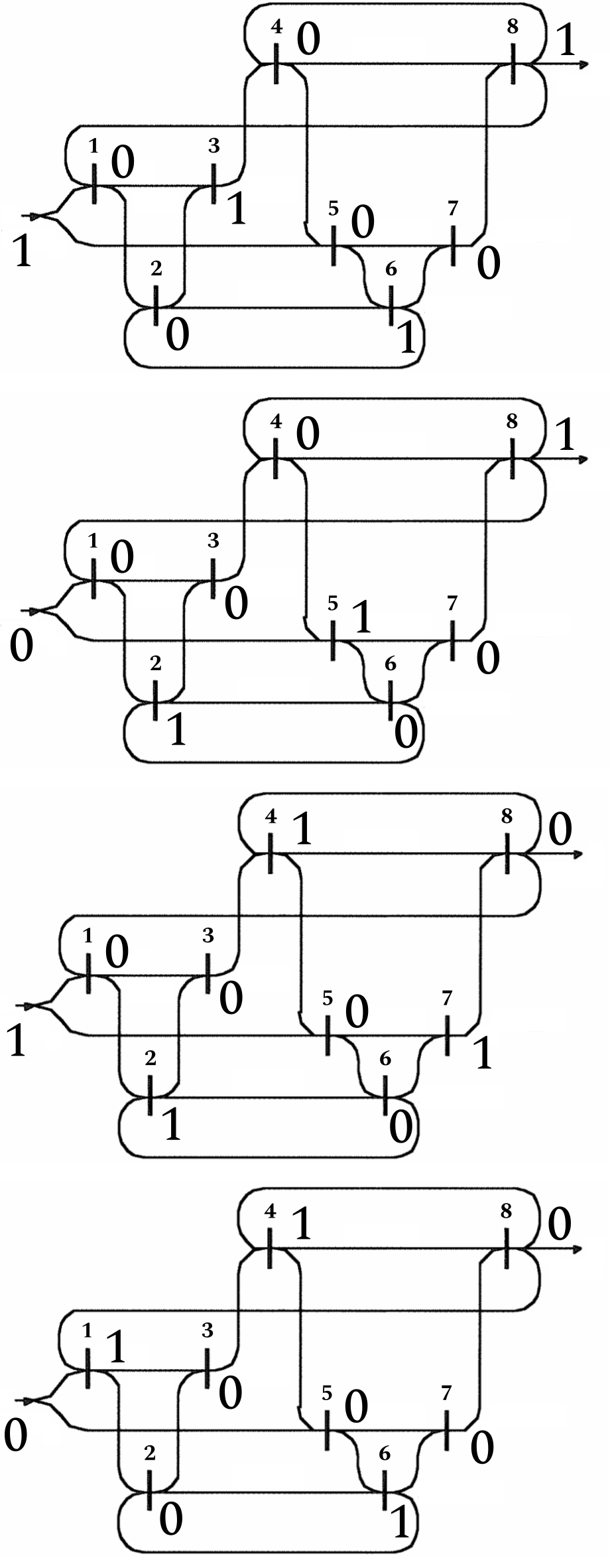}
        \end{figure}

\newpage


\noindent The waveforms that emerge globally at each marker are as follows:

\begin{figure}[h!]
            \centering
            \includegraphics[scale=0.18]{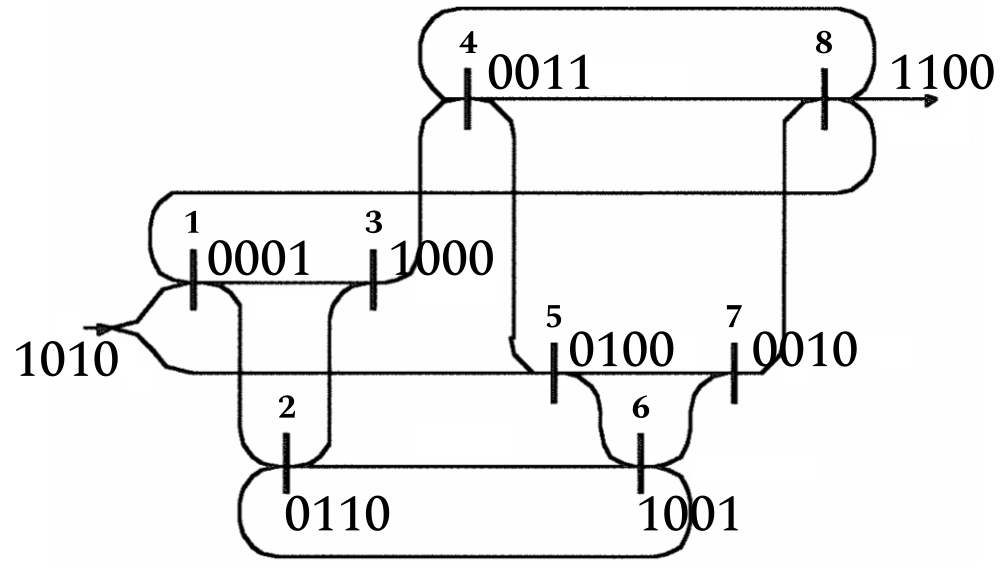}
        \end{figure}

\noindent The following table shows both ways of reading the modulator. Reading the table row-wise shows what happens ``locally" at each marker at a given balanced state, and reading the table column-wise shows the waveform that emerges ``globally" at a given marker. The transition details are not shown in this table.

\begin{center}
\begin{tabular}{c c c c c c c c c c c c c}
    Input & & $\mathbf{1}$ & $\mathbf{2}$ & $\mathbf{3}$ & $\mathbf{4}$ & $\mathbf{5}$ & $\mathbf{6}$ & $\mathbf{7}$ & $\mathbf{8}$ & & Output\\ 
    $1$ & & $0$ & $0$ & $1$ & $0$ & $0$ & $1$ & $0$ & $1$ & & $1$\\
    $0$ & & $0$ & $1$ & $0$ & $0$ & $1$ & $0$ & $0$ & $1$ & & $1$\\ 
    $1$ & & $0$ & $1$ & $0$ & $1$ & $0$ & $0$ & $1$ & $0$ & & $0$\\
    $0$ & & $1$ & $0$ & $0$ & $1$ & $0$ & $1$ & $0$ & $0$ & & $0$
\end{tabular}
\end{center}


\noindent We now redraw E4 in a form similar to that conceived by Jeff James.\footnote{Here, we reference excerpts from ``Spencer-Brown's Counter: A simple analysis by Jeff James", a document that is now archived on the internet. Note that Jeff James drew a somewhat different circuit than the one shown above.} Markers labeled 1-8 are circles labeled 1-8 in this graph. Readers are encouraged to check that this graph is isomorphic to E4:





\begin{figure}[h!]
            \centering
            \includegraphics[scale=0.15]{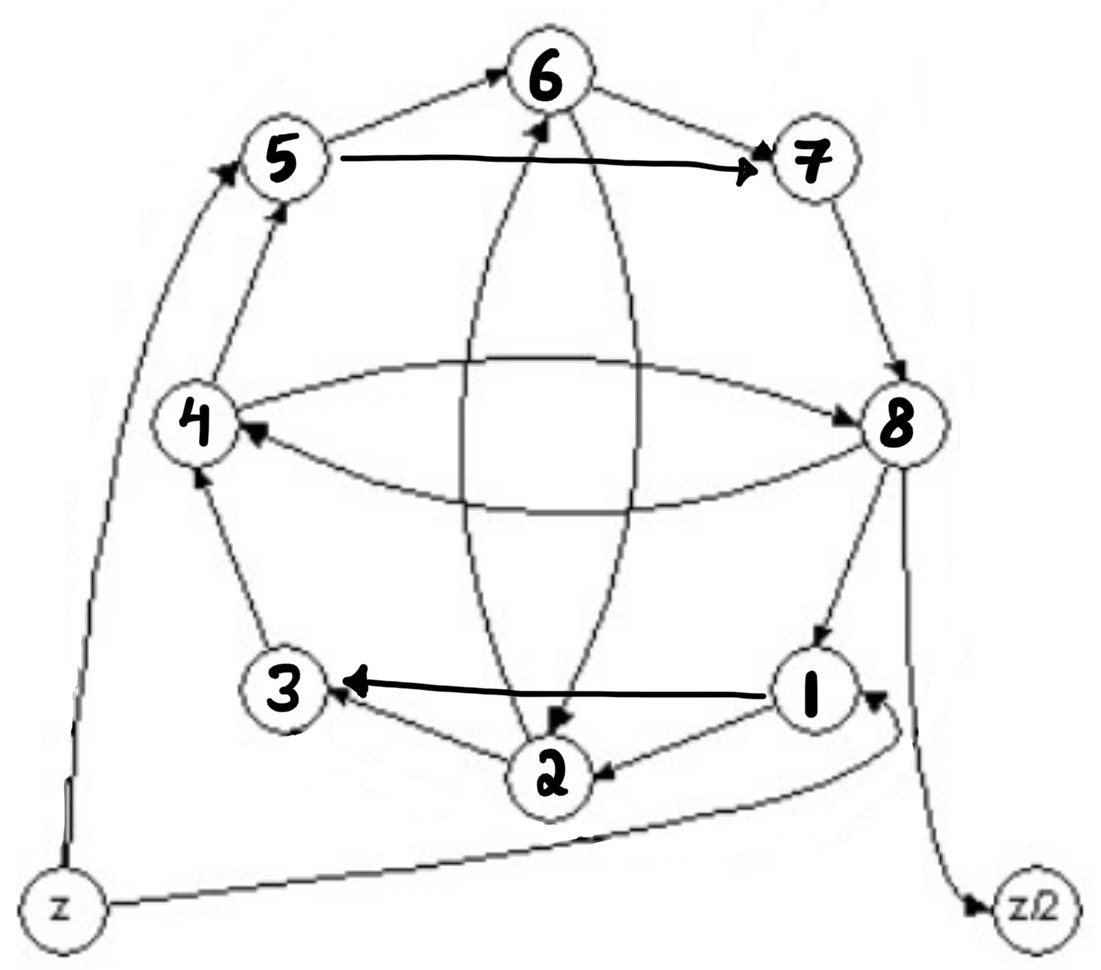}
        \end{figure}


\noindent While E4 in this form is a sight to marvel at, the local and global ways of looking at E4 in this form is even more astonishing, revealing all its symmetries.\footnote{In the words of Jeff James: ``The four phases of the counter can be seen below to rotate a triangular constellation of black nodes (\textit{marked 1}) around the circle. The states progress clockwise, as does the constellation within them. Each state represents a settling of the circuit given a new input."} First we see what happens locally at each marker at each successive balanced state:


%

\begin{figure}[h!]
            \centering
            \includegraphics[scale=0.44]{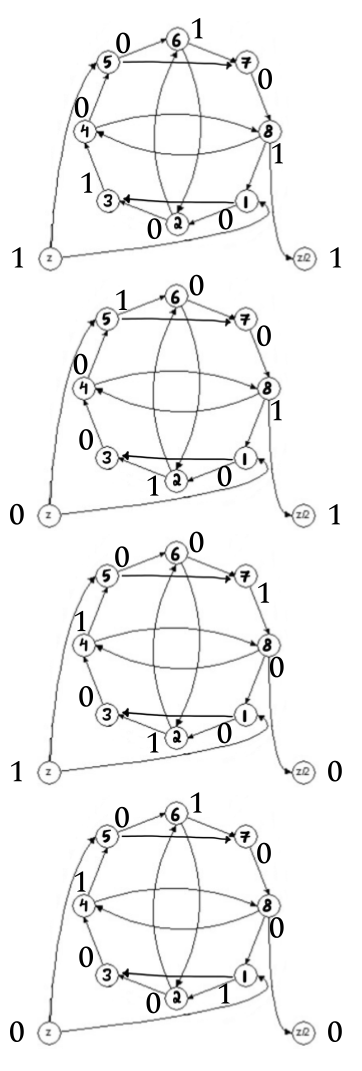}
        \end{figure}

\newpage

\noindent Next, we see the waveforms that emerge globally at each marker:

\begin{figure}[h!]
            \centering
            \includegraphics[scale=0.4]{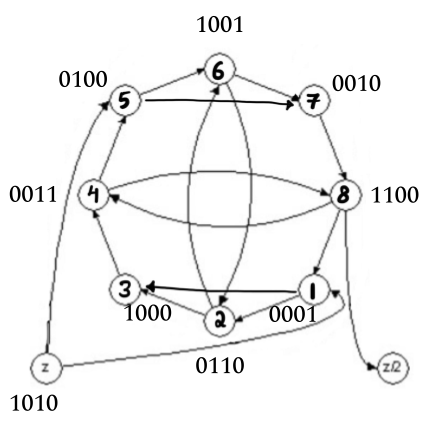}
        \end{figure}







\subsection{Binary sieving with E4}\hphantom{.}

\noindent We now consider the binary forms of the labels of the markers in E4 and the form of the periodic binary sequences that emerge at these labels. Observe the Jeff James form of E4. First note that the prime markers have single count bitstrings and composite markers have composite unit bitstrings. Consider even markers $2, 4, 6, 8$, and note markers $2$ and $6$, and markers $4$ and $8$ are in opposition. Their binary waveforms (marker 4 = 0011, marker 2 = 0110, marker 6 = 1001, marker 8 = 1100) correspond to numbers $3$, $6$, $9$, and $12$, multiples of $3$. Now consider markers $1, 3, 5, 7$. Their binary waveforms (marker 1 = 0001, marker 7 = 0010, marker 5 = 0100, marker 3 = 1000) correspond to $1, 2, 4, 8$, powers of $2$. Note, the number $10$ corresponds to the input waveform $1010$, so the holes created by E4 are the primes $5,7, 11$. Markers 1 and 5 (corresponding to squares 1 and 4) are in opposition, creating 2 successive spaces that yield the primes $2$ and $3$. These primes are used as strikers to obtain all primes up to $25=5^2$ (Lemma 1A).

\begin{center}
\begin{tabular}{c c c c c c | c c c c c c c}
    Input & & $\mathbf{1}$ & $\mathbf{2}$ & $\mathbf{3}$ & $\mathbf{4}$ & $\mathbf{5}$ & $\mathbf{6}$ & $\mathbf{7}$ & $\mathbf{8}$ & & Output\\ 
    $1$ & & $0$ & $0$ & $1$ & $0$ & $0$ & $1$ & $0$ & $1$ & & $1$\\
    $0$ & & $0$ & $1$ & $0$ & $0$ & $1$ & $0$ & $0$ & $1$ & & $1$\\
    \hline
    $1$ & & $0$ & $1$ & $0$ & $1$ & $0$ & $0$ & $1$ & $0$ & & $0$\\
    $0$ & & $1$ & $0$ & $0$ & $1$ & $0$ & $1$ & $0$ & $0$ & & $0$
\end{tabular}
\end{center}

\noindent Observe again, the table that illustrates the local (row-wise) and global (column-wise) ways of reading the modulator E4. Reading the modulator E4 row-wise, we see a gap between markers 4 and 5, splitting E4 into two distinct sections. Also note the gap between the second and the third stable state (second row to third row). Observe how opposing quadrants are repeated exactly. Reflexivity balances both these half sections together. Indeed, E4 self-regulates by its own structure. 

\subsection{Modulator structure in general and binary counting} \hphantom{.}

\noindent Note that the general idea of a modulator is that it can accept a waveform or varying input and it will output a waveform of half of the input frequency. Thus if $...10101010...$ is the input then the output would be $...11001100...$. To see how this works, imagine a black box $\mathbf{M}$ that takes an input $x$ and gives an output $y$. Write $y\mathbf{M}x$ for this situation and call it a \textit{stable state} of $\mathbf{M}$. It is assumed that once $\mathbf{M}$ receives $x$, there will be a time interval after the output $y$ is stably computed. Here are the rules for $\mathbf{M}$:

\begin{itemize}
    \item If $y\mathbf{M}0$ is a stable state of $\mathbf{M}$ and $0$ changes to $1$, then $y\mathbf{M}1$ is a stable state of $\mathbf{M}$. In other words, if the input $0$ changes to $1$, no change happens to the output.
    \item If $y\mathbf{M}1$ is a stable state of $\mathbf{M}$ and $1$ changes to $0$, then $y'\mathbf{M}1$ is the resulting stable state, where $0'=1$ and $1'=0$. In other words, if the input $1$ changes to $0$, then the output changes (either from $0$ to $1$ or from $1$ to $0$).
    \item One can start the modulator either in any of the states $0\mathbf{M}0$, $0\mathbf{M}1$,$1\mathbf{M}0$,$1\mathbf{M}1$.
\end{itemize}

\noindent Let's see how $\mathbf{M}$ behaves. Start $\mathbf{M}$ in the state $0\mathbf{M}0$ and change the input successively, following the rules. We find:
\begin{center}
\noindent $0\mathbf{M}0$\\
$0\mathbf{M}1$\\
$1\mathbf{M}0$\\
$1\mathbf{M}1$\\
$0\mathbf{M}0$\\
$...$
\end{center}

\noindent As you can see, $\mathbf{M}$ does indeed behave as advertised. If we put a sequence of modulators together, then the pattern of the outputs corresponds to binary counting:

\begin{center}
\noindent $0\mathbf{M}0\mathbf{M}0\mathbf{M}0$\\
$0\mathbf{M}0\mathbf{M}0\mathbf{M}1$\\
$0\mathbf{M}0\mathbf{M}1\mathbf{M}0$\\
$0\mathbf{M}0\mathbf{M}1\mathbf{M}1$\\
$0\mathbf{M}1\mathbf{M}0\mathbf{M}0$\\
$0\mathbf{M}1\mathbf{M}0\mathbf{M}1$\\
$0\mathbf{M}1\mathbf{M}1\mathbf{M}0$\\
$0\mathbf{M}1\mathbf{M}1\mathbf{M}1$\\
$1\mathbf{M}0\mathbf{M}0\mathbf{M}0$\\
$1\mathbf{M}0\mathbf{M}0\mathbf{M}1$\\
$1\mathbf{M}0\mathbf{M}1\mathbf{M}0$\\
$1\mathbf{M}0\mathbf{M}1\mathbf{M}1$\\
$1\mathbf{M}1\mathbf{M}0\mathbf{M}0$\\
$1\mathbf{M}1\mathbf{M}0\mathbf{M}1$\\
$1\mathbf{M}1\mathbf{M}1\mathbf{M}0$\\
$1\mathbf{M}1\mathbf{M}1\mathbf{M}1$
\end{center}

\noindent Note that at the next step we go back to $0\mathbf{M}0\mathbf{M}0\mathbf{M}0$. This concatenation of three modulators counts to $16$. We see that from this, and the fact that the modulator $\mathbf{M}$ can be built from operations of distinction in the circular pattern of the designs of Chapter 11 (so called ``JK-Flip Flop circuits" in the engineering parlance), numerical arithmetic arises from recursive acts of distinction.

\subsection{The Spencer-Brown sieve} \hphantom{.}

\noindent It is the purpose of this section to point out the general structure of modulators and to point out that this structure is intimately related with the design and behavior of the Spencer-Brown sieve (the GSB sieve). We have discussed this version of the Sieve of Eratosthenes while constructing the ``sect pyramid" from Lemma 3 earlier in the paper. Here we start again from the beginning.

\noindent Write two place holders $**$.
There are two of them. Add one to get $3$ and place $3$ on either side of the two stars.
$$3**3$$
\noindent Now count that find $4$ symbols. Add $1$ and get $5$. Place $5$.
$$53**35$$
\noindent Count gives $6$ and add one gives $7$. Place $7$. Note that we are generating prime numbers.
$$753**357$$
\noindent But now the $3$ can skip to the left by $3$ or to the right by $3$ and occupy an outer place. When this can be done it must be done.
$$3753**3573$$
\noindent One can try to skip $5$ by $5$ steps, but it will go too far, not at the exact ends. So we are done and we can Count and add one to get 10+1=11, the next prime number.
$$11\hphantom{.}3753**3573\hphantom{.}11$$
\noindent Neither $3$ nor $5$ can skip to the ends and so we can count and add one to get $12+1 = 13$.
$$13\hphantom{.}11\hphantom{.}3753**3573\hphantom{.}11\hphantom{.}13$$
\noindent Now $3$ or $5$ can skip forward and strike the ends. When more than one strike can occur, we indicate the smallest one.
$$3\hphantom{.}13\hphantom{.}11\hphantom{.}3753**3573\hphantom{.}11\hphantom{.}13\hphantom{.}3$$
\noindent At this point no further strikes are available and so we can Count and add one.
$$17\hphantom{.}3\hphantom{.}13\hphantom{.}11\hphantom{.}3753**3573\hphantom{.}11\hphantom{.}13\hphantom{.}3\hphantom{.}17$$
\noindent This process continues forever, producing prime numbers. We can think of the GSB sieve as a generalized modulator whose stable states occur when no more strikes are available. At such points the output is the next prime in form of Count $+ 1$. Then the modulator can be stablilized again.

\newpage

\noindent Here is what we did in one chart.
                 $$**$$
                $$3**3$$
               $$53**35$$
             $$3753**3573$$
         $$11\hphantom{.}3753**3573\hphantom{.}11$$
    $$3\hphantom{.}13\hphantom{.}11\hphantom{.}3753**3573\hphantom{.}11\hphantom{.}13\hphantom{.}3$$
$$17\hphantom{.}3\hphantom{.}13\hphantom{.}11\hphantom{.}3753**3573\hphantom{.}11\hphantom{.}13\hphantom{.}3\hphantom{.}17$$
$$...$$

\noindent The next chart shows how the original modulator patter is inside this GSB Sieve:
                 $$*\#$$
                $$3*\#3$$
               $$53\#*35$$
             $$3753\#*3573$$
         $$11\hphantom{.}3753*\#3573\hphantom{.}11$$
    $$3\hphantom{.}13\hphantom{.}11\hphantom{.}3753*\#3573\hphantom{.}11\hphantom{.}13\hphantom{.}3$$
$$17\hphantom{.}3\hphantom{.}13\hphantom{.}11\hphantom{.}3753\#*3573\hphantom{.}11\hphantom{.}13\hphantom{.}3\hphantom{.}17$$
$$...$$

\noindent One can see a modulator pattern in the center stars. These can be regarded as filled with primes repeating in modulator pattern ready to come out to the boundaries after the strikes have come to stability. 

\newpage

\noindent In the chart below we fill in those central primes. It is on the second appearance of each central prime that it goes outward to take its place in the external row of the GSB Sieve.
                    $$\LoF{\vphantom{@}}{\LoF{\vphantom{@}}{5}}\hphantom{.}\LoF{\vphantom{@}}{3}$$
                $$3  \LoF{\vphantom{@}}{5}\LoF{\vphantom{@}}{\LoF{\vphantom{@}}{7}} \hphantom{.}      3$$
               $$53  \LoF{\vphantom{@}}{\LoF{\vphantom{@}}{11}}\hphantom{.}\LoF{\vphantom{@}}{7} \hphantom{.}     35$$
             $$3753 \LoF{\vphantom{@}}{11}\LoF{\vphantom{@}}{\LoF{\vphantom{@}}{13}} \hphantom{.}   3573$$
         $$11\hphantom{.}3753 \LoF{\vphantom{@}}{\LoF{\vphantom{@}}{17}}\hphantom{.}\LoF{\vphantom{@}}{13} \hphantom{.}  3573\hphantom{.}11$$
    $$3\hphantom{.}13\hphantom{.}11\hphantom{.}3753 \LoF{\vphantom{@}}{17}\LoF{\vphantom{@}}{\LoF{\vphantom{@}}{19}} \hphantom{.}  3573\hphantom{.}11\hphantom{.}13\hphantom{.}3$$
$$17\hphantom{.}3\hphantom{.}13\hphantom{.}11\hphantom{.}3753 \LoF{\vphantom{@}}{\LoF{\vphantom{@}}{23}}\hphantom{.}\LoF{\vphantom{@}}{19} \hphantom{.}  3573\hphantom{.}11\hphantom{.}13\hphantom{.}3\hphantom{.}17$$

\noindent The Spencer-Brown sieve (the GSB sieve) is a powerful evocation of the generation of the prime numbers. Looking at it in this way one can see that it is part of the intricate story of how numbers are generated by circuits of distinctions.

\noindent The modulation structure in the central multisectors emerge from two distinct placeholders, i.e., ``two successive spaces of expression" (Chapter 11, \textit{Laws Of Form}).
$$\LoF{\vphantom{@}}{*}\hphantom{.}*$$
$$\LoF{\LoF{\vphantom{@}}{\hphantom{.}}}{\!}\hphantom{.}\LoF{\vphantom{@}}{\hphantom{.}}$$

\fbox{\begin{minipage}{32em}
\noindent Every fact about every prime in the universe is contained in our knowledge of the situation of the first prime $2$. Furthermore, we do not even need to know that $2$ is a prime, merely that it is the first multiplicative integer. $1$ is nonexistent as a multiplier, divisor, exponent, or root. Thus $2$ can tell us the answer to every question that can possibly be asked about primes.
\end{minipage}}

\noindent Following this revelatory excerpt from page 191 of Appendix 8, we close by restating Spencer-Brown's more-difficult theorem, which he ``first proved in the summer of 1998".

\fbox{\begin{minipage}{32em}
\noindent \textbf{Theorem 2}:\\
For all natural squares, there is a prime in every half segment.
\end{minipage}}

\newpage

\chapter*{Glossary}

\noindent \textbf{square segment, or \textit{segment}}: the set of whole numbers between successive squares.\\
\textbf{seg$\mathbf{(n)}$} - the set of whole numbers between $n^2$ and $(n+1)^2$, where $n \in \mathbb{N}$.\\
\noindent \textbf{seg$\mathbf{(x)}$:} the set of whole numbers between $x^2$ and $(x+1)^2$, where $x\in\mathbb{R}$.\\
$^{\circ}\pi$\textbf{seg}$\mathbf{(x)}$: the number of primes between $x^2$ and $(x+1)^2$.


\noindent \textbf{Theorem 1:} There is at least one odd prime in the square segment of every natural number.

\noindent \textbf{Lemma 1:} A number $h$ is prime \textit{iff} $h$ is an integer and $h$ does not make an integer quotient with any integer $i$ in case $2\leq i \leq \sqrt{h}$.

\noindent \textbf{pg of $\mathbf{h}$}, or \textbf{prime generators of the primes in $\mathbf{h}$:} the set of primes $\leq \sqrt{h}$.

\noindent \textbf{arithmetic $\mathbf{A}$:} the system of integers\\
\textbf{system $\mathbf{S}$:} the mapping of the system of integers onto a line of regularly spaced indistinguishable marks or points stretching endlessly in either direction.

\noindent \textbf{prime multisectors $\mathbf{d_j}$:} operators on the system $S$ that imbed with or strike every second, every third, every fifth, ... every $p_j$\textsuperscript{th} element, where $p_j$ is the $j$\textsuperscript{th} prime number. In other words, a prime multisector is all multiples of a given prime. \\
\textbf{decimation/decimator:} synonymous with multisection/multisector.

\noindent $\mathbf{A'}$: the odd integers.\\
$\mathbf{S'}$: the subsystem of $S$ that maps to the odd integers.

\noindent \textbf{strike:} (verb) In operating the sieve, we choose a prime number and ``strike" out all multiples of that prime number. The operation of the imbedding of a multisector onto the elements of system $S$: when we speak of imbedding a multisector, this corresponds to striking the corresponding set of points in $S$. (noun) A representative piece of $S$ that results from the striking of a set of multisectors. (collective noun) The primorial or prime factorial, also called hash or factorial strike. The primorial of $n$ is the product of primes $\leq n$ and the factorial strike is the strike using these primes. This is also called the full reflexive strike.

\noindent $\mathbf{p_{\sim j}}$\textbf{:} the set of all primes $\leq p_j$.\\
$\mathbf{d_{\sim j}}$\textbf{:} the set of all prime multisectors $\leq d_j$.

\noindent \textbf{Figure 1:} The \textbf{full reflexive strike} (or the \textbf{prime paradigm}) of the odd primes $3, 5$. The ringed points are mirror-points at multiples of $2$, which do not strictly belong to the odd-strike.

\noindent \textbf{quantum or half-strike:} half of the primorial used to generate the prime paradigm.


\noindent \textbf{Lemma 2:} There are exactly $n$ odd numbers between $n^2$ and $(n+1)^2$.


\noindent \textbf{sect:} a consecutive sequence of composite numbers.\\
\noindent \textbf{R-sect or reflexive sect:} a ``palindromic" odd sect constructed using Recipe 1.\\
\noindent \textbf{Q-sect or quadrantic sect:} an odd sect constructed using Recipe 2.

\noindent \textbf{Recipe 1:} A recipe to create reflexive sects or R-sects. Beginning at a reference mark ($0$) in a system $S$, this recipe creates a palindromic odd sect by imbedding or striking every $p_j$\textsuperscript{th} element (where $p_j$ is the $j$\textsuperscript{th} prime number) on either side of the reference, including two central points corresponding to $-1$ and $+1$.\\
\textbf{Recipe 2:} A recipe to create quadrantic or Q-sects. Beginning at a reference mark ($0$) in a system $S$, this recipe creates an odd sect by imbedding or striking every $p_j$\textsuperscript{th} element (where $p_j$ is the $j$\textsuperscript{th} prime number) on the right of the reference (zero), including one point corresponding to $+1$.

\noindent \textbf{mR$\mathbf{(d_{\sim j})}$:} the maximum R-sect constructed using prime multisectors $d_1, d_2, ... d_j$.\\
\textbf{mQ$\mathbf{(d_{\sim j})}$:} the maximum Q-sect constructed using prime multisectors $d_1, d_2, ... d_j$.

\noindent \textbf{Lemma 3:} The maximum possible sect for the pg of every $n^2$ is mR$(d_{\sim j}) = d_{j-1} -1$ or mQ$(d_{\sim j}) = \frac{d_{j+1}-1}{2}-1$. An exception to Lemma 3 is $n = 23$.



\noindent \textbf{Lemma 4:} mR$(d_{\sim j}) \geq$ mQ$(d_{\sim j})$. In other words, $p_{j-1}-1\geq \frac{p_{j+1}-1}{2}-1$, which when rearranged implies $p_{j+1}\leq 2p_{j-1}+1$.


\noindent \textbf{paradigmatic place (pp) and nonparadigmatic place (np):} a multisector is in a paradigmatic place, or pp, if it imbeds in a point that can be identified with the associated prime in $A'$. Otherwise it is in a nonparadigmatic place np. 



\noindent \textbf{Theorem 1A:} If $x$ is a real number $\geq 1$, the next prime $> x$ is at a maximum distance of $2\sqrt{x}$ from $x$.


\noindent \textbf{Table 1:} A fake prime table illustrating maximum distances between fake primes corresponding to Lemma 4 being ``barely true", i.e., $p_{j+1} = 2p_{j-1}+1$. 

\noindent \textbf{Table 2:} A fake prime table that compares actual maximum R-sects (determined by actual primes listed in column 1) with fake maximum Q-sects for fake successive primes that are barely within the limits of Theorem 1A (where a fake successive prime is the nearest odd number within two square-roots of the actual prime in column 1).



\noindent \textbf{elemental:} A set of consequential theorems that dovetail with each other in such a way that they stand or fall together.


\noindent \textbf{Proposition (5):} $p_{i+1} \leq p_i + \frac{1}{k}\sqrt{p_i} \quad \quad (k = \frac{1}{2})$.


\noindent \textbf{Ratchet point:} A point on the number line where $k$ in Proposition (5) advances to a new value from which it cannot slip back.

\noindent \textbf{Theorem 1B:} For any arbitrarily large $k$, Proposition (5) is true for all numbers from some $i$ upwards.

\noindent \textbf{Corollary:} If $i \geq 32, (p_{i+1} - p_i)$ and $(p_i - p_{i-1})$ are both less than $|p_i^{\frac{1}{2}}|$

\noindent \textbf{Theorem 2:} For all natural squares, there is a prime in every half segment.




\newpage

\chapter*{Appendix 1:\\ Sharper Conjectures on Primes Between Squares}

\noindent Spencer-Brown conjectures that the number of primes between $n^2$ and $(n+1)^2$ is approximately $\frac{n}{\log n}$. This appendix provides sharper conjectures on the existence of primes between squares as provided by Spencer-Brown.

\noindent In Appendix 7, he writes: ``The total number $t(n)$ of primes between $n^2$ and $(n+1)^2$, where $n$ is a natural number greater than $1$, is within the limits $A - (B - 1) < t(n) < A + (B - 1)$ where $A$ is $\frac{n}{\log n}$ and $B$ is $\frac{A}{\log A}$. Since my formula, which I call the Prime Limit Theorem, is not necessarily true unless $n$ is one of the natural numbers $2, 3, 4, ...,$ it further confirms the existence of a fundamental relationship between logarithms to the base $e$ and the natural numbers, their squares, and the primes" (pages 182-183, \textit{Laws Of Form}, Bohmeier Verlag, 2015).

\noindent ``I failed to notice in Appendix 7 of Reference 1 that my upper limit to the possible number of primes in seg$(n)$, \hphantom{..} $A+(B-1) \quad (A=\frac{n}{\log n}, B=\frac{A}{\log A})$, is also an upper limit to $\pi(n)$" (page 200, \textit{Laws Of Form}, Bohmeier Verlag, 2015).

\noindent In \textit{Notes} to Appendix 8, (page 201), Spencer-Brown writes: ``The best upper limit to $\pi(x)$ hitherto was Sylvester's 1892 refinement of Tchebychef's result, notably that $0.95695 \frac{x}{\log x} < \pi(x) < 1.04423 \frac{x}{\log x}$ for all sufficiently large $x$. We may ignore the lower limit since it is easy to prove that $\frac{n}{\log n} < \pi(n)$ for all $n \geq 17$. Sylvester's upper limit is too small until $n$ is of the order of $10^{11}$, whereas my upper limit is true for all $n$" (pages 201-202, \textit{Laws Of Form}, Bohmeier Verlag, 2015). 

\noindent The following excerpts are taken from the handwritten version of \textit{Primes Between Squares}, reproduced as follows:

\noindent ``More importantly my lower limit $A - (B - 1) < t$ where $A$ is $\frac{n}{\log n}$, $B$ is $\frac{A}{\log A}$, and $t$ is the total number of primes between $n^2$ and $(n+1)^2$, provides a much more realistic answer for the minimum number of primes between $n$ and $2n$, and even this can be improved if we make accurate assessments that avoid the approximations due to using logarithms. My upper limit $A+(B-1)>\pi(n)$ is so near to the actual count that a simple modification of it will give an approximation comparable to Riemann's for up to $10^7$.

\noindent My approximation (one of several I composed in 1998) is $Qx = a + b - c = \pi(x)$ where $a = \frac{x}{\log x}, b = \frac{a}{\log a}, c = \frac{jb}{\log jb}$ and $j$ is an adjustor between $0$ and $1$. Setting $j = \frac{1}{2}$ gives a very good estimate for $\pi(x)$ in the early reaches of $x$. My more sophisticated formula (2) in Reference 1 will match Riemann's results as far as we wish to take it, but the summation using a hand calculator is time-consuming for large $x$. I have since improved it by invoking the subroutine SB$(x)=$ Li $(x)+$ Li $((\sqrt{x}-1)^2)/2 - 1.5$. When $x$ is large we can substitute SB$([\sqrt{x}+t]^2)$ for the summation term in (2), calling the result SB$_t(x)$. The difference to the outcome will be negligible. Thus S$_{1.32}(10^{11}) = 4118052495.21$ and SB$_{1.32}(10^{11}) = 4118052495.46$ show hardly any difference, except that the former took 12 hours to compute, and the latter only 3 minutes. The Riemann formula for the same $x$ yields R$(10^{11}) = 4118052494.46$ and is extremely laborious to compute. All the results in Reference 1 were achieved with a programmable hand calculator, using the simplest and most accurate formula for Li$(x):$ Li$(x) = \gamma + \log\log x + \sum_{k=1}^{\infty}\frac{\log^k x}{k! k}$. 

\noindent No one hitherto seems to have remarked that all such approximations to $\pi(x)$ are meaningless unless we state over what stretch of numbers they are claimed to be reasonably accurate. Riemann's and my approximations are claimed to be reasonably accurate near zero, but Li$(x)$ would of course provide a much closer approximation to $\pi(x)$ in the region of a Littlewood crossing\footnote{See Appendix 9, \textit{Laws of Form} (Bohmeier Verlag, 2015).}. In fact the Riemann-Ramanujan formula\footnote{Original footnote: The disadvantage of the Riemann-Ramanujan formula is that it requires what is theoretically an infinite number of computations of logarithmic integrals, whereas my subroutine above, whose results converge to those of the R-R formula, requires only two such computations to achieve what is practically the same answer.} is merely a complicated modification of Li$(x)$ to make it a closer approximation to $\pi(x)$ in the region of zero".

\chapter*{Appendix 2: A Constellation of Conjectures}

\noindent Here, we list a number of well-known conjectures. We will not provide references, since the readers can see these conjectures in a variety of number theory expositions.

\noindent \textbf{Bertrand's postulate}\footnote{Although it is called a postulate, it is a theorem that has been proved by Tchebychef and Ramanujan.}: For any integer $n>3$, there exists at least one prime number $p$ such that $n < p < 2n-2$.

\noindent \textbf{Legendre's conjecture:} For every positive integer $n$, there is a prime number between $n^2$ and $(n+1)^2$.

\noindent \textbf{Oppermann's conjecture:} For every positive integer $n$, there is a prime between $n(n-1)$ and $n^2$, and a prime between $n^2$ and $n(n+1)$.

\noindent \textbf{Brocard's conjecture:}  For all $n$, if $p_n$ is the $n$\textsuperscript{th} prime number, there are at least four prime numbers between $p_{n}^2$ and $p_{n+1}^2$.

\noindent \textbf{Andrica's conjecture:} For all $n$, if $p_n$ is the $n$\textsuperscript{th} prime number, $\sqrt{p_{n+1}}-\sqrt{p_n}>1$.

\noindent \textbf{Cramer's conjecture:} For all $n$, if $p_n$ is the $n$\textsuperscript{th} prime number, $p_{n+1}-p_n = O((\log p_n)^2)$.


The following three conjectures are also on prime gaps $g_n = p_{n+1}-p_n$:
\begin{figure}[h!]
            \centering
            \includegraphics[scale=0.53]{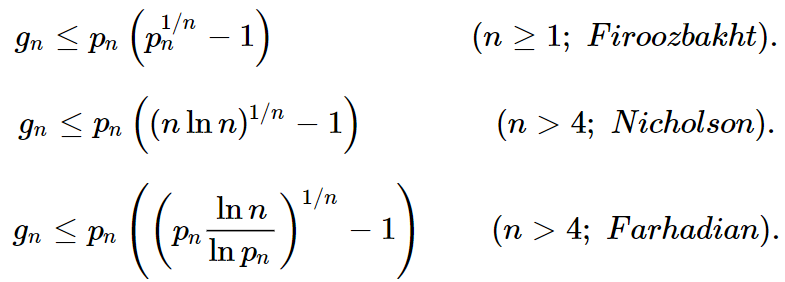}
        \end{figure}

\noindent \textbf{Forgues' conjecture}: For all $n$, if $p_n$ is the $n$\textsuperscript{th} prime number, $\left(\frac{\log(p_{n+1})}{\log(p_n)}\right)^n < e$.


\noindent The form of Appendix 8 connects Bertrand's postulate, as well as Legendre’s, Oppermann’s, Brocard’s, Andrica's, Firoozbakht's, Cramer's, Farhadian's, Nicholson's, and Forgues' conjectures. 

\noindent The ratchet points in Theorem 1A and Theorem 1B appear to follow ratchet points in Andrica’s conjecture OEIS A084974 ``Primes that show the slow decrease in the larger values of the Andrica function $Af(k) = \sqrt{p_{k+1}} - \sqrt{p_k}$, where $p_k$ denotes the $k$’th prime". Also, see\footnote{OEIS A182514 ``Primes $prime(n)$ such that $\frac{prime(n+1)}{prime(n)}^n > n$".\\ OEIS A383591 ``Smallest prime $p$ where the absolute difference of the gaps to the adjacent primes exceeds $n\log(p)$".\\ OEIS A111943 ``Prime $p$ with prime gap $q - p$ of $n$-th record Cramer-Shanks-Granville ratio, where $q$ is smallest prime larger than $p$ and C-S-G ratio is $\frac{(q-p)}{(log p)^2}$".\\ OEIS A079098 ``Conjectured values of greatest $k$ such that for any consecutive primes $q, q'$, $k\leq q<q', \sqrt{q'}-\sqrt{q}<\frac{1}{n}$".} OEIS A182514, A383591, A111943, and A079098.


\noindent Denote:\\
Oppermann's conjecture = $O$\\
The standard Andrica conjecture = $a$ or Theorem 1A (Spencer-Brown)\\
The strong Legendre conjecture = $L$\\
The standard Legendre conjecture = $l$ or Theorem 1 (Spencer-Brown).\\
The strong Brocard conjecture = $B$\\
The standard Brocard conjecture = $b$\\
Here $x \to y$ means $x$ implies $y$.

\noindent Germ\'{a}n Andr\'{e}s Paz (April 2014) constructs Conjecture 1 such that:\\ 
Conjecture 1 $\to l$\\
Conjecture 1 $\to L$\\
Conjecture 1 $\to b$\\
Conjecture 1 $\to O$\\
Conjecture 1 $\to a$.\\
Paz (June 2014) also finds reentrant reformulations of Conjecture 1.

\noindent Matt Visser (2019) constructs a stronger form of the standard Andrica conjecture $A$ and shows $A \to O$. In ``Verifying other weaker conjectures", Visser shows $O \to L$, $L \to l$, $L \to B$, $B \to b$, $O \to a$, $a \to l$, and $l \to$ the weak Andrica conjecture. 

\begin{figure}[h!]
            \centering
            \includegraphics[scale=0.45]{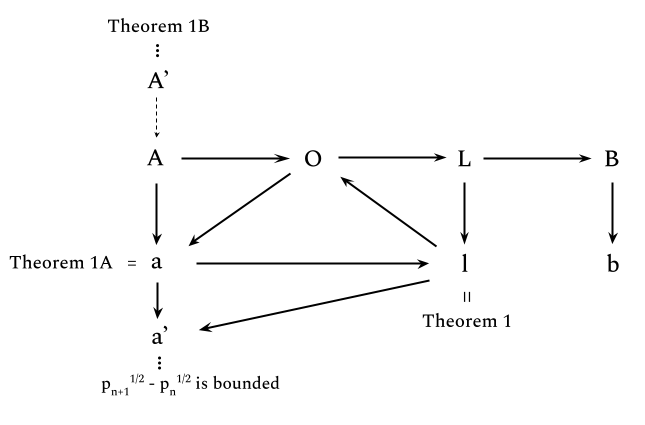}
        \end{figure}

\noindent Theorem 1 (Spencer-Brown) is the standard Legendre Conjecture $l$.\\
Theorem 1A (Spencer-Brown) is the standard Andrica conjecture $a$.\\

\begin{figure}[h!]
            \centering
            \includegraphics[scale=0.45]{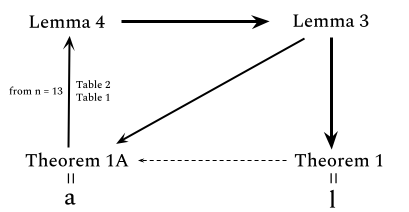}
        \end{figure}

\noindent Theorem 1B (Spencer-Brown) pushes the strengthening of Andrica's conjecture to its limit. The following statements are under investigation:
Theorem 1B and its Corollary (Spencer-Brown) imply $A$ and Conjecture 1.\\ 
Theorem 2 (Spencer-Brown) is a strengthening of $O$.

\noindent It is fascinating to observe how Theorem 1B, the Corollary, and Theorem 2, relate to the constellation of conjectures on prime gaps.

\chapter*{Appendix 3: Exception at $23$}

\noindent Spencer-Brown's Lemma 3 states that the maximum odd sect is either the R-sect mR$(d_{\sim j})=d_{j-1}-1$ or Q-sect mQ$(d_{\sim j})=\frac{(d_{j+1}-1)}{2}-1$ with a single exception at $23$.

\noindent For odd primes, $d_{\sim 9}=$ the set of primes $d_2$ through $d_9$ $ = \{3,5,7,9,13,17,19,23\}$.\\
Hence, mR$(d_{\sim 9})=19-1=18$, and mQ$(d_{\sim 9})=\frac{(29-1)}{2}-1=13$.

\noindent The maximum odd R-sect using $d_{\sim 9}$, written mR$(d_{\sim 9})$, is of the form:
\begin{figure}[h!]
            \centering
            \includegraphics[scale=0.6]{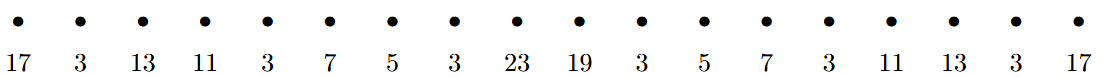}
        \end{figure}

\noindent The maximum odd Q-sect using $d_{\sim 9}$, written mQ$(d_{\sim 9})$, is of the form:
\begin{figure}[h!]
            \centering
            \includegraphics[scale=0.76]{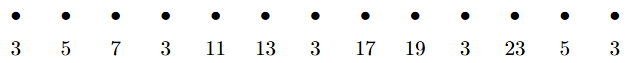}
        \end{figure}

\noindent However, neither of these is the largest odd-sect using $d_{\sim 9}=\{3,5,7,11,13,17,19,23\}$. On page 193, Appendix 8, \textit{Laws of Form}, Bohmeier Verlag, Spencer-Brown writes: ``In fact Lemma 3 identifies the maximum possible sect for the pg of every $n^2$ except in case $n=23$, for which (uniquely) there are six maximum odd sect of $19$. ..."

\noindent On each of the following two pages, we provide seven maximum odd sects and their mirror images that involve primes $3,5,7,11,13,17,19,23$. One of these is the maximum odd R-sect mR$(d_{\sim 9})$ of length $18$ as shown above. This can be identified by a dashed line that splits the sect in the middle. Note, this dashed line is located at a multiple of $17\#$. The other six sects on each of the following two pages are the anomalous maximum odd sects of length $19$, one greater than the length stipulated by Lemma 3.



\newpage
\thispagestyle{empty}

\begin{landscape}

\begin{figure}[h!]
            \centering
            \includegraphics[scale=0.83]{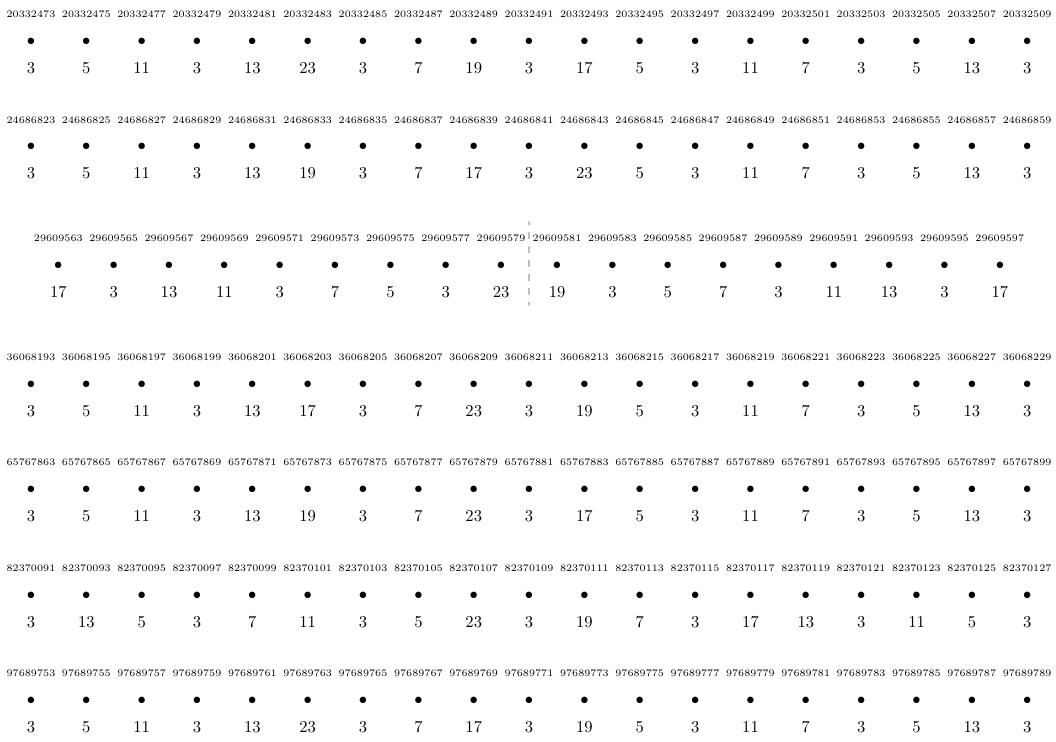}
        \end{figure}

\end{landscape}

\newpage
\thispagestyle{empty}

\begin{landscape}

\begin{figure}[h!]
            \centering
            \includegraphics[scale=0.75]{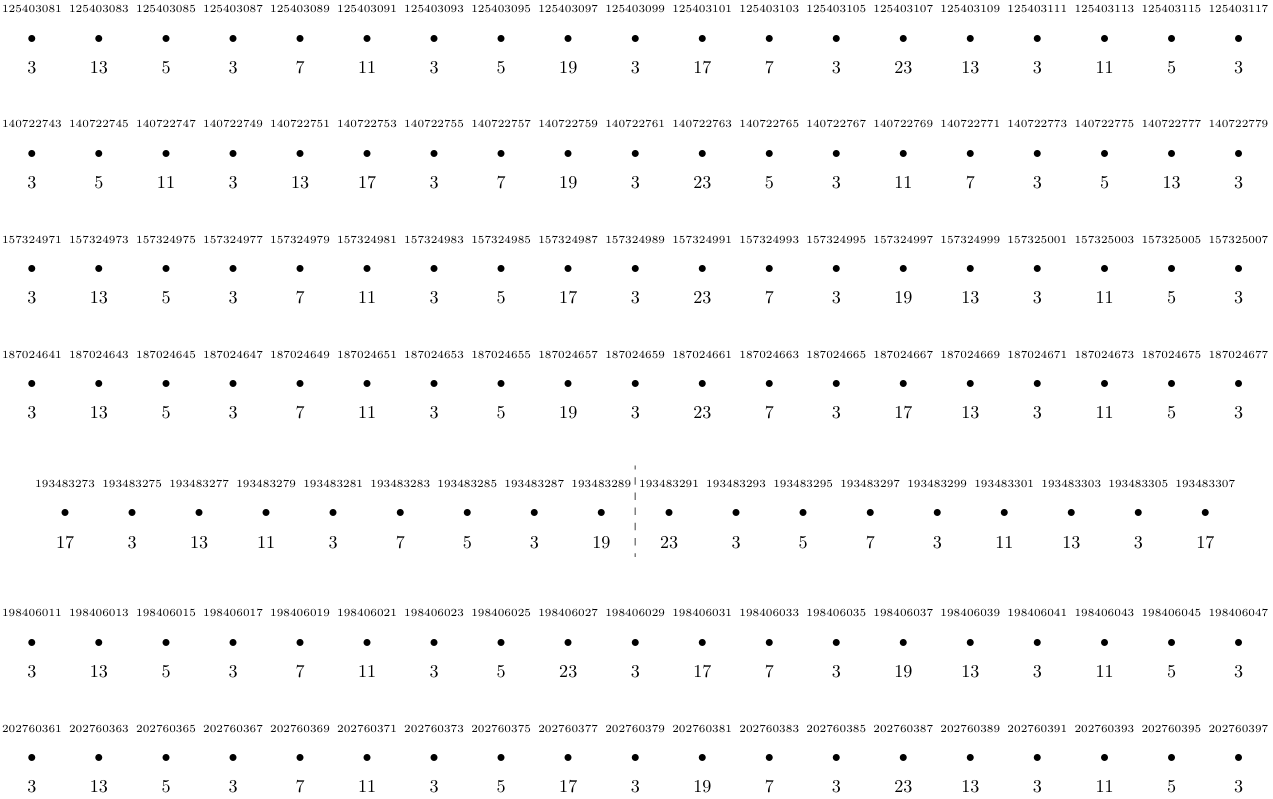}
        \end{figure}

\end{landscape}

\newpage


\noindent On page 201 of \textit{Laws Of Form} (Bohmeier Verlag, 2015), Spencer-Brown writes: ``My theorem that Lemma 3 always denotes a maximum sect has this single exception at $d_{\sim 9}$. The fact never ceases to astonish me, because it means that the primes $\leq 23$ can be rearranged to strike more efficiently than they do in the paradigm. It would seem \textit{prima facie} to be impossible, but the fact that it can happen once makes it even more surprising that it can never happen again".\footnote{The anomaly at $23$ deserves a paper all on its own (Original footnote, pg. 201, \textit{Laws Of Form}, 2015).}

\noindent Spencer-Brown also points out the anomaly at $23$ in an inequality for Euler's totient function $\phi(n)$, the number of positive integers up to $n$ that are relatively prime to $n$. This inequality, provided by J.B. Rosser and L. Schoenfeld in ``Approximate formulas for some functions of prime numbers" (1962), is as follows:
$$\phi(n)\geq e^{-\gamma}\frac{n}{\log\log n + \frac{5}{2 e^\gamma \log\log n}}$$
for all $n\geq 3$, with the single exception of $n=23\#$, for which the constant $\frac{5}{2}$ must be replaced by $2.50637$. Here $\gamma$ is Euler's constant\footnote{Euler's constant (or the Euler-Mascheroni constant) $\gamma = \int_{1}^{\infty}\left(-\frac{1}{x}+\frac{1}{\floor{x}}\right)dx\approx 0.57721...$}.

\noindent Yet another reference to the anomaly at $23$ is in the study of finite field extensions of the rational numbers, which is ``indispensable to the number theory, even if one's ultimate goal is to understand properties of diophantine expressions and equations in the ordinary integers" (Avigad, 2004). In the Translator's Introduction to ``Dedekind's 1871 version of the theory of ideals", Avigad writes: ``It can happen, however, that the ``integers" in such extensions fail to satisfy unique factorization, a property that is central to reasoning about the ordinary integers. In 1844, Ernst Kummer observed that unique factorization fails for the cyclotomic integers with exponent $23$, i.e., the ring $\mathbb{Z}[\zeta]$ of integers in the field $\mathbb{Q}(\zeta)$, where $\zeta$ is a primitive twenty-third root of unity" (Avigad, 2014). \footnote{\href{https://www.andrew.cmu.edu/user/avigad/Papers/ideals71.pdf}{https://www.andrew.cmu.edu/user/avigad/Papers/ideals71.pdf}}


\noindent That such anomalies exist at $23$ is worth the special attention of number theorists.

\newpage

\noindent Conceived decades before ``Primes Between Squares", the 19\textsuperscript{th} and 23\textsuperscript{rd} poems in \textit{23 Degrees of Paradise} written by Spencer-Brown and published under the pseudonym James Keys, anticipate and conjure the deep mathematical insight and mystery of 23.

\noindent The 19\textsuperscript{th} poem in \textit{23 Degrees of Paradise} is titled ``23". The 23\textsuperscript{rd} poem is titled ``PS":


\noindent ``\textit{P S}

\noindent An idle poet sees and writes\\
A book of wonders and delights\\
So full of beauty, life, and love\\
That poetry's what he writes it of.

\noindent Next comes a busy college don\\
Who looks into the poet's delights\\
And then a book of learning writes,\\
And poetry's what he writes it on.

\noindent And finally a student, who\\
Has read the don and learned to scoff,\\
Cries `Is this all that poets do?'\\
Takes his degree and writes it off.

\noindent The years pass. A little child,\\
Rude, untutored, undefiled,\\
Idle, full of joy, and wild,\\
Comes to where the Poet Smiled,

\noindent And from the vision of it writes\\
A book of wonders and delights\\
So full of beauty, life, and love\\
That poetry's what he writes it of."\\
(James Keys, \textit{23 Degrees of Paradise})


\chapter*{References:}

\noindent Andrica, D. ``Note on a conjecture in prime number theory". \textit{Studia Univ. Babe\c{s}-Bolyai Math.} \textbf{31} (4): 44-48.

\noindent Avigad, J. ``Dedekind's 1871 Version of the Theory of Ideals". Carnegie Mellon University Repository. \href{https://philpapers.org/rec/AVIDV}{https://philpapers.org/rec/AVIDV}. March 19, 2004.

\noindent Baker, R. C. Harman, G. Pintz, J. ``The Difference Between Consecutive Primes, II". Proceedings of the London Mathematical Society, 83: 532-562. 2001. \\ https://doi.org/10.1112/plms/83.3.532

\noindent Ellis, N. \textit{Awakening Osiris}: A New Translation of the Egyptian Book of the Dead. Phanes Press. 1988.

\noindent Copi, I. M. \textit{The Theory of Logical Types}. Routledge \& Kegan Paul. London. 1971.


\noindent Dickson, L. E. \emph{History of the Theory of Numbers}. I. Divisibility and Primality. 1952.

\noindent Feliksiak, J. ``The Nicholson's Conjecture". ScienceOpen Preprints. 24 February 2021. 

\noindent Flagg, J. M. Kauffman, L. H. Sahoo, D. ``Laws Of Form and the Riemann Hypothesis". \textit{Laws of Form - A Fiftieth Anniversary}. World Scientific. 2021.

\noindent James, J. ``Spencer-Brown's Counter". \href{https://web.archive.org/web/20100306005801/http://www.lawsofform.org/counter/index.html}{https://web.archive.org/web/20100306005801/\linebreak http://www.lawsofform.org/counter/index.html}

\noindent Kauffman, L. H. ``Network Synthesis and Varela's Calculus". Int. J. General Systems, Vol. 4, pp. 179-187. 1977.

\noindent Kauffman, L. H. ``Modulators and Imaginary Values". \textit{Laws of Form - A Fiftieth Anniversary}. World Scientific. 2021.

\noindent Keys, J. \emph{23 Degrees of Paradise}. W Heffer \& Sons Ltd. Cambridge. 23 April, 1970.

\noindent Keys, J. \textit{Only Two Can Play This Game}. The Julian Press, Inc. 1972.

\noindent Landau, Edmund. ``Gel\"{o}ste und ungel\"{o}ste Probleme aus der Theorie der \linebreak Primzahlverteilung und der Riemannschen Zetafunktion". Jahresbericht der Deutschen Mathematiker-Vereinigung 21: 208-228. \href{http://eudml.org/doc/145337}{http://eudml.org/doc/145337}. 1912.

\noindent Legendre, A. M. \textit{Essai sur la Th\'{e}orie des Nombres}. Seconde \'{E}dition. Quatri\`{e}me Partie. IX. ``D\'{e}monstration de divers th\'{e}or\`{e}ms sur les progressions arithm\'{e}tiques". Paris, Chez Courcier, Imprimeur-Libraire pour les Math\'{e}matiques, quai des Augustins, n\textsuperscript{o} 57. 1808.

\noindent Nagura, J. ``On the interval containing at least one prime number". Proc. Japan Acad. 28, no. 4, 177–181. 1952.



\noindent Online Encyclopedia of Integer Sequences (OEIS). A084974 ``Primes that show the slow decrease in the larger values of the Andrica function $Af(k)=\sqrt{p_{k+1}}-\sqrt{p_k}$, where $p_k$ denotes the $k$-th prime".

\noindent Panaitopol, L. ``Intervals Containing Prime Numbers". NNTDM 7, 4, 111-114. 2001.

\noindent Paz, G. A. ``On Legendre's, Brocard's, Andrica's, and Oppermann's Conjectures". \linebreak \href{https://arxiv.org/abs/1310.1323}{https://arxiv.org/abs/1310.1323}. 2 April 2014.

\noindent Paz, G. A. ``Two statements that are equivalent to a conjecture related to the \linebreak distribution of prime numbers". \href{https://arxiv.org/abs/1406.4801}{https://arxiv.org/abs/1406.4801}. 19 June 2014.

\noindent Petzold, C. \textit{The Annotated Turing: A Guided Tour through Alan Turing's Historic Paper on Computability and the Turing Machine}. Wiley Publishing, Inc. 2008.


\noindent Rosser, J. B. Schoenfeld, L. ``Approximate formulas for some functions of prime \linebreak numbers". Illinois J. Math, 6, 64-94. 1962.

\noindent Smith, H. ``Andrica's conjecture". Retrieved from the Wayback Machine. \linebreak \href{https://www.oocities.org/hjsmithh/PrimeSR/index.html}{https://www.oocities.org/hjsmithh/PrimeSR/index.html}.

\noindent Spencer-Brown, G. ``Commentary on Appendix 9". Unpublished manuscript. 2007.

\noindent Spencer-Brown, G. \emph{Gesetze der Form}. Bohmeier Verlag. 1997.

\noindent Spencer-Brown, G. ``Haploid and Diploid Numbers". Unpublished manuscript. 1998.

\noindent Spencer-Brown, G. ``Introduction to Reductors". Unpublished manuscript. December 1992.

\noindent Spencer-Brown, G. \emph{Laws Of Form}. Revised Seventh English edition. Bohmeier Verlag. 2020.

\noindent Spencer-Brown, G. \emph{Laws Of Form}. Revised Sixth English edition. Bohmeier Verlag. 2015.

\noindent Spencer-Brown, G. ``Primes Between Squares". Handwritten edition. Copy 22. June 2000.

\noindent Spencer-Brown, G. ``The Universal Operators in Logic". Unpublished manuscript.

\noindent Stein, M. L. Ulam, S. M. Wells, M. B. ``A Visual Display of Some Properties of the Distribution of Primes". Amer. Math. Monthly 71, 516-520. 1964

\noindent Visser, M. ``Strong Version of Andrica's Conjecture". International Mathematical Forum, Vol. 14, no. 4, 181-188. Hikari Ltd. \href{https://arxiv.org/abs/1812.02762}{https://arxiv.org/abs/1812.02762}. 2019.








\end{document}